\documentclass[]{article}

\usepackage{amsmath, amsthm, amssymb, mathtools}
	\mathtoolsset{showonlyrefs}
\usepackage[a4paper, margin=2.54cm]{geometry}
\usepackage{parskip} 
\usepackage{booktabs} 
\usepackage{ifthen}
\usepackage{url}
\usepackage{hyperref}
	\hypersetup{colorlinks, linkcolor={green!70!black}, citecolor={lightred},urlcolor={azure}} 
\usepackage{algorithm}
\usepackage{algorithmicx}
\usepackage{algpseudocode}
	\algnewcommand\algorithmicforeach{\textbf{for each}}
	\algdef{S}[FOR]{ForEach}[1]{\algorithmicforeach\ #1\ \algorithmicdo}
	\algnewcommand{\IfThen}[2]{\State \algorithmicif\ #1\ \algorithmicthen\ #2}
	\algnewcommand{\LineComment}[1]{\State \textcolor{azure}{\texttt{/*\;#1\;*/}}}
	\algrenewcommand\algorithmiccomment[1]{\hfill \textcolor{azure}{\texttt{//\;#1}}}
	\algrenewcommand{\textproc}{\texttt}
	\newlength{\WidthOfInput}
	\newlength{\WidthOfOutput}
\usepackage{colortbl} 
\usepackage{tcolorbox} 
	\tcbuselibrary{most}
\usepackage{caption} 
	\DeclareCaptionLabelFormat{cont}{#1 #2 (continued)}
\usepackage{tikz}
	\usetikzlibrary{calc, arrows.meta, positioning, patterns, external, decorations.pathmorphing, decorations.pathreplacing, decorations.markings, shadings, shapes.arrows, shapes.geometric, matrix, fadings, backgrounds}
	\tikzexternaldisable
\usepackage{pgfplots}
	\pgfplotsset{compat=newest}
	\usepgfplotslibrary{groupplots, fillbetween}
	\pgfkeys{/pgf/number format/.cd,1000 sep={\,}} 
	\pgfplotscreateplotcyclelist{color list 2}{lightred, lightblue, orange, lightgreen, oker, brown, teal, cyan, magenta, violet, gray}
	\newlength{\figurewidth}
	\newlength{\figureheight}
\usepackage{pgfplotstable}
\usepackage{multirow}
\usepackage{bbm} 
\usepackage[labelformat=simple]{subcaption}
	
\usepackage{mathrsfs}
\usepackage{setspace}
\definecolor{azure}{rgb}{0.0, 0.5, 1.0}
\definecolor{olive}{RGB}{0,151,76}
\definecolor{oker}{RGB}{234,170,45}
\definecolor{lightred}{RGB}{237,20,35}
\definecolor{lightblue}{RGB}{45,121,181}
\definecolor{bordeaux}{RGB}{127,57,86}

\colorlet{lightgreen}{green!70!black}
\colorlet{yescolor}{lightgreen}
\colorlet{nocolor}{lightred}

\tikzset{
	default line/.style={line cap=round, smooth, thick, mark=*, mark options={solid}, mark size=.8pt},
	hollow line/.style={default line, mark options={fill=white}},
	dotted line/.style={default line, dashed, dash pattern=on 0pt off 4\pgflinewidth},
	time step/.style={thick, line cap=round},
	kinetic/.style={thick, line cap=round, lightred},
	diffusive/.style={thick, line cap=round, lightblue, dashed, dash pattern=on 0pt off 4\pgflinewidth},
	collision/.style={circle, fill=lightred, inner sep=1pt},
	>={Stealth[width=1.5mm,length=1.5mm]},
	default arrow/.style={->, semithick, line cap=round},
	bar legend entry/.style={mark=*, only marks, mark size=2},
	large dot/.style={circle, inner sep=3pt},
	yes dot/.style={large dot, fill=yescolor},
	no dot/.style={large dot, fill=nocolor},
}

\pgfplotsset{%
	default axis/.style={%
		width=\figurewidth,
		height=\figureheight,
		major tick length=2pt,
		minor tick length=2pt,
		every tick/.style={black, line cap=round},
		tick label style={font=\scriptsize},
		axis on top,
		legend style={draw=none, font=\scriptsize, at={(0.97,0.03)}, anchor=south east, fill=none, legend cell align=left},
	},
	loglog axis/.style={%
		default axis,
		cycle list name=color list 2,
		xmode=log,
		ymode=log,
	},
	samples axis/.style={%
		default axis,
		every tick/.style={opacity=0},
		every y tick label/.style={xshift=10pt},
		xlabel=$j$,
		ybar=-9pt,
		bar width=9pt,
		ylabel={$N_{\level_j}$},
		legend style={draw=none, font=\scriptsize, at={(1.1,1.0)}, anchor=north east, fill=none, legend cell align=left},
		bar legend,
		ymajorgrids,
		major grid style={white, line width=.75pt},
		x axis line style={opacity=0},
		y axis line style={opacity=0},
	},
	/pgfplots/bar legend/.style={%
    		/pgfplots/legend image code/.code={%
     			\draw[##1,/tikz/.cd,bar width=.5em] plot coordinates { (11pt,0em)};}
	},
}
\newtcolorbox[auto counter, number within=section]{todobox}[1][]{
	outer arc=2pt,
	arc=1pt,
	colframe=red,
	colback=red!20,
	attach boxed title to top left={yshift=-3mm, xshift=3mm},
	enhanced,
	boxed title style={
		colback=red,
		top=3pt,
		bottom=3pt,
	},
  	fonttitle=\sffamily,
	top=12pt,
  	title=TO DO~\thetcbcounter,
}

\newtcolorbox[auto counter, number within=section]{infobox}[1][]{
	outer arc=2pt,
	arc=1pt,
	colframe=olive,
	colback=olive!20,
	attach boxed title to top left={yshift=-3mm, xshift=3mm},
	enhanced,
	boxed title style={
		colback=olive,
		top=3pt,
		bottom=3pt,
	},
  	fonttitle=\sffamily,
	top=12pt,
  	title=INFO~\thetcbcounter,
}

\newtcolorbox[auto counter, number within=section]{notebox}[1][]{
	outer arc=2pt,
	arc=1pt,
	colframe=azure,
	colback=azure!20,
	attach boxed title to top left={yshift=-3mm, xshift=3mm},
	enhanced,
	boxed title style={
		colback=azure,
		top=3pt,
		bottom=3pt,
	},
  	fonttitle=\sffamily,
	top=12pt,
  	title=NOTE~\thetcbcounter,
}

\newcommand{\slopetriangle}[4]{
	\pgfplotsextra{ %
		\pgfkeysgetvalue{/pgfplots/xmin}{\xmin}
		\pgfkeysgetvalue{/pgfplots/xmax}{\xmax}
		\pgfkeysgetvalue{/pgfplots/ymin}{\ymin}
		\pgfkeysgetvalue{/pgfplots/ymax}{\ymax}

		\pgfmathsetmacro{\xArel}{#1}
		\pgfmathsetmacro{\yArel}{#3}
		\pgfmathsetmacro{\xBrel}{#1-#2}
		\pgfmathsetmacro{\yBrel}{\yArel}
		\pgfmathsetmacro{\xCrel}{#4<0?\xBrel:\xArel}

		\pgfmathsetmacro{\lnxB}{\xmin*(1-(#1-#2))+\xmax*(#1-#2)}
		\pgfmathsetmacro{\lnxA}{\xmin*(1-#1)+\xmax*#1}
		\pgfmathsetmacro{\lnyA}{\ymin*(1-#3)+\ymax*#3}
		\pgfmathsetmacro{\lnyC}{\lnyA+abs(#4)*(\lnxA-\lnxB)}
		\pgfmathsetmacro{\yCrel}{\lnyC-\ymin)/(\ymax-\ymin)}
        
		\coordinate (A) at (rel axis cs:\xArel,\yArel);
		\coordinate (B) at (rel axis cs:\xBrel,\yBrel);
		\coordinate (C) at (rel axis cs:\xCrel,\yCrel);

		\pgfmathsetmacro{\factor}{#4<0?0.6:0.1}
		\pgfmathtruncatemacro{\absslope}{abs(#4)}
		\draw[]   (A) --node[pos=\factor,yshift=1ex,xshift=-0.5ex] {\scriptsize \absslope} (B) -- (C) -- cycle;
    }
}

\newcommand{\figref}[1]{Figure~\ref{#1}}
\newcommand{\tabref}[1]{Table~\ref{#1}}

\newcommand{\secref}[1]{Section~\ref{#1}}

\renewcommand{\algref}[1]{Algorithm~\ref{#1}}

\newcommand{\kdstep}[3]{
	\draw[kinetic] (#1, 0) -- ({#1+#2}, 0);
	\draw[diffusive] ({#1+#2}, 0) -- ({#1+#3}, 0);
	\node[collision] at ({#1+#2}, 0) {};
} 

\newcommand{\bsy}{{\boldsymbol{y}}}

\newcommand{\bbE}{{\mathbb{E}}}

\newcommand{\bbV}{{\mathbb{V}}}

\DeclareSymbolFont{bbold}{U}{bbold}{m}{n}
\DeclareSymbolFontAlphabet{\mathbbold}{bbold}

\newcommand{\calC}{{\mathcal{C}}}

\newcommand{\calE}{{\mathcal{E}}}

\newcommand{\calL}{{\mathcal{L}}}

\newcommand{\calN}{{\mathcal{N}}}
\newcommand{\calO}{{\mathcal{O}}}

\newcommand{\calQ}{{\mathcal{Q}}}

\providecommand{\argmin}{\operatorname*{arg\,min}}

\newcommand{\VAtimestepno}{i}

\newcommand{\VAeventno}{k}
\newcommand{\VAtimestep}{\delta t}
\newcommand{\VAkintime}{\tau}
\newcommand{\VAdiftime}{\theta}

\newcommand{\VAaddeventno}[2]{#1_#2}
\newcommand{\VArnvel}{\nu}
\newcommand{\VArndif}{\chi}
\newcommand{\VAfactor}{M}
\newcommand{\VArnexp}{\epsilon}
\newcommand{\VArate}{R}

\newcommand{\VAbackgroundspeedmean}{\mu_v}
\newcommand{\VAbackgroundspeedstdv}{\sigma_v}

\newcommand{\VAfineindex}[1]{\kappa_{#1}}

\newcommand{\E}[1]{\bbE\mathopen{}\left[{#1}\right]\mathclose{}}
\newcommand{\V}[1]{\bbV\mathopen{}\left[{#1}\right]\mathclose{}}
\newcommand{\cost}[1]{\calC\mathopen{}\left({#1}\right)\mathclose{}}

\newcommand{\mse}[1]{\text{MSE}\mathopen{}\left({#1}\right)\mathclose{}}
\newcommand{\order}[1]{\calO\mathopen{}\left({#1}\right)\mathclose{}}

\newcommand{\level}{\ell}

\newcommand{\VAendtime}{T}
\newcommand{\VAparticleno}{n}
\newcommand{\VAncollisions}{K}
\newcommand{\VArandparam}{\bsy}
\newcommand{\VAqoi}{Q}
\newcommand{\VAqoiest}{\calQ}

\newcommand{\B}[1]{\mathsf{B_#1}}

\newcommand{\kiessymbool}{\zeta} 

\title{Multilevel asymptotic-preserving Monte Carlo for kinetic-diffusive particle simulations of the Boltzmann-BGK equation}
\author{Bert Mortier\footnote{%
KU Leuven, Department of Computer Science, NUMA Section.\newline
\protect\makebox[1.8em]{}Celestijnenlaan 200A box 2402, 3001 Leuven, Belgium.\newline
\protect\makebox[1.8em]{}\href{mailto:bert.mortier@kuleuven.be}{\{bert.mortier, pieterjan.robbe, giovanni.samaey\}@kuleuven.be}
}
\and Pieterjan Robbe\footnotemark[1]{}
\and Martine Baelmans\footnote{%
KU Leuven, Department of Mechanical Engineering, Applied Mechanics and Energy Conversion Section.\newline
\protect\makebox[1.8em]{}Celestijnenlaan 300A box 2421, 3001 Leuven, Belgium.\newline
\protect\makebox[1.8em]{}\href{mailto:tine.baelmans@kuleuven.be}{tine.baelmans@kuleuven.be}
}
\and Giovanni Samaey\footnotemark[1]{}
}

\date{}

\begin{document}

\maketitle

\begin{abstract}We develop a novel multilevel asymptotic-preserving Monte Carlo method, called Multilevel Kinetic-Diffusion Monte Carlo (ML-KDMC), for simulating the kinetic Boltzmann transport equation with a Bhatnagar--Gross--Krook (BGK) collision operator. This equation occurs, for instance, in mathematical models of the neutral particles in the plasma edge of nuclear fusion reactors. In this context, the Kinetic-Diffusion Monte Carlo method is known to maintain accuracy both in the low-collisional and the high-collisional limit, without an exploding simulation cost in the latter. We show that, by situating this method within a Multilevel Monte Carlo (MLMC) framework, using a hierarchy of larger time step sizes, the simulation cost is reduced even further. The different levels in our ML-KDMC method are connected via a new and improved recipe for correlating particle trajectories with different time step sizes. Furthermore, a new and more general level selection strategy is presented. We illustrate the efficiency of our ML-KDMC method by applying it to a one-dimensional test case with nonhomogeneous and anisotropic plasma background. Our method yields significant speedups compared to the single-level KDMC scheme, both in the low and high collisional regime. In the high-collisional case, our ML-KDMC outperforms the single-level KDMC method by several orders of magnitude.
\end{abstract}



\section{Introduction}\label{sec:intro}

Kinetic equations play a vital role in many modern applications. For example, in mathematical models for nuclear fusion reactors such as ITER and DEMO, see~\cite{ITER}, the physics of neutral particles in the plasma is modeled using the Boltzmann transport equation with a Bhatnagar--Gross--Krook (BGK) collision operator, see, e.g.,~\cite{bhatnagar54, lux91, stangeby00}. The kinetic equation then boils down to simulating every individual collision of the neutral particle with the plasma background. Near the plasma edge, there is an area of increased neutral-plasma collision rates. This reduces the heat load on the the plasma-facing components significantly. However, the increase in the number of collisions has severe implications on the computational burden of the kinetic description. On the other hand, it is well-known that, in the high-collisional limit, the behanviour of the neutral particles converges to an advection-diffusion process, see, e.g.,~\cite{othmer00}. The latter can be simulated cheaply using biased random walks.

Thus, there are regions in the domain where a kinetic description is required, and also regions with a high collision rate where this kinetic description becomes intractable, but where a diffusive approximation exists that is cheap to simulate. In many works, domain decomposition is the method of choice for solving these type of problems, see, e.g.~\cite{boyd11, densmore07}. However, the domain decomposition approach requires a good partitioning of the domain into a kinetic and a diffusive part, and an efficient coupling between both. Alternative hybrid approaches, that avoid this coupling altogether, are the so-called \emph{asymptotic-preserving Monte Carlo} (APMC) methods~\cite{pareschi01, dimarco18}. These methods use a single approximation scheme throughout the domain, such that the method has the accuracy of the kinetic simulation in the low-collision regions, and the efficiency of a diffusive simulation in the high-collision regions. Asymptotic-preserving methods were originally developed in the context of radiation transport, see, e.g.,~\cite{fleck71, fleck84}, and later on also for neutron transport, see~\cite{borgers92}, and the Boltzmann-BGK equation, see, e.g.,~\cite{gabetta97, crouseilles04, crestetto12, crestetto18, radtke13}.

One example of such an asymptotic-preserving scheme can be found in~\cite{dimarco18}. This method uses an implicit time discretization to obtain an unconditionally stable fixed time step Monte Carlo method, thereby limiting the simulation cost in the high-collisional limit. However, this method has a computational cost that grows unboundedly with decreasing time step size. This problem is alleviated in the Kinetic-Diffusion Monte Carlo (KDMC) method, an asymptotic-preserving Monte Carlo method developed in~\cite{mortier20}. The latter method uses hybridized particles that exhibit both kinetic behanviour and diffusive behanviour depending on the local collisionality.

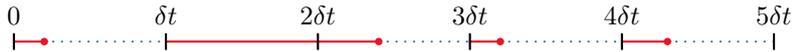
\begin{figure}
\centering
\begin{tikzpicture}[xscale=2]
	\foreach \i in {0, 1, ..., 5}{
		\draw[time step] (\i, -.1) -- (\i, .1);
		\node[anchor=south, outer sep=3pt] at (\i, 0) {$\ifthenelse{\i=0}{0}{\ifthenelse{\i=1}{}{\i}{\VAtimestep}}$};
	}
	\begin{scope}[on background layer]
		\kdstep{0}{.2}{1}
		\kdstep{1}{1.4}{2}
		\kdstep{3}{.2}{1}
		\kdstep{4}{.3}{1}
	\end{scope}
\end{tikzpicture}
\caption{An illustration of a particle trajectory in the KD scheme. A kinetic step (\raisebox{2pt}{\protect\tikz{\protect\draw[kinetic] (0,0) -- (.5,0);}}) ends with a collision (\raisebox{1pt}{\protect\tikz{\protect\node[collision] {};}}), and is always followed by a diffusive step (\raisebox{2pt}{\protect\tikz{\protect\draw[diffusive] (0,0) -- (.5,0);}}) until the beginning of the next time step. Not that, for the particle trajectory shown above, no collision takes place in $[\VAtimestep, 2\VAtimestep)$.}
\label{fig:ml_sim_timedomain}
\end{figure}

In this paper, we present a multilevel extension of the KDMC scheme. This extension is termed \emph{Multilevel Kinetic-Diffusion Monte Carlo} (ML-KDMC). Our method combines the asymptotic-preserving kinetic-diffusion scheme from~\cite{mortier20} with the \emph{Multilevel Monte Carlo} (MLMC) method, see, e.g.,~\cite{giles08, giles15}. MLMC methods use a hierarchy of coarse approximations to reduce the computational cost of a simulation. In the context of particle simulations, the hierarchy of coarse approximations can be constructed by subsequently increasing the time step size in the simulation. The goal of our ML-KDMC method is then to reduce the computational cost of the kinetic simulation, while keeping the flexibility of the KDMC method.

Earlier work on combining asymptotic-preserving particle methods with multilevel Monte Carlo methods can be found in~\cite{loevbak19}. Our current work differs from~\cite{loevbak19} in several ways. First, in~\cite{loevbak19}, the APMC scheme of~\cite{dimarco18} is used, whereas the present work uses the KDMC scheme from~\cite{mortier20}. Second, our multilevel method uses a new recipe to generate correlated particle trajectories with different time step sizes, a notoriously difficult problem that was identified in~\cite{loevbak20}. This is a consequence of the choice for the KDMC scheme from~\cite{mortier20}. Third, we apply our method to the relevant case of a non-homogeneous plasma background, and present a new level selection strategy. Finally, we show that our method can be applied without change to an anisotropic plasma background. We remark that the ideas presented in this paper can also be combined with the APMC scheme of~\cite{dimarco18}, see~\cite{loevbak20a}.

The remainder of this text is organized as follows. First, in~\secref{sec:ml_sim}, we briefly discuss the KDMC scheme from~\cite{mortier20}. Next, in~\secref{sec:ml_cor}, we present the main contribution of this work, i.e., the improved recipe for correlating particle trajectories. In~\secref{sec:mlapmc}, we discuss the MLMC method in the context of kinetic equations, and address the construction of an optimal hierarchy of coarse approximations for multilevel sampling. Finally, in~\secref{sec:numres}, we present numerical results that illustrate the superiority of our ML-KDMC scheme over the standard KDMC method from~\cite{mortier20} in terms of computational cost.

\section{The KD simulation scheme\label{sec:ml_sim}}

In this section, we discuss the kinetic-diffusion (KD) simulation scheme for the Boltzmann-BGK equation introduced in~\cite{mortier20}. We first present the general idea. Suppose the time domain of the simulation is discretized into disjoint time intervals of equal length $\VAtimestep$. Let $x(t)$ and $v(t)$ denote the position and velocity of a particle at time $t$. In the KD simulation scheme, particles alternatingly follow kinetic and diffusive trajectories, as shown in Figure~\ref{fig:ml_sim_timedomain}. Particles move kinetically with a constant velocity $v(t)$, until a collision occurs. The kinetic step is oblivious to any time discretization, meaning that a kinetic trajectory may span several time steps, as illustrated in Figure~\ref{fig:ml_sim_timedomain}. Then, for the remainder of the time step in which the collision occurred, the particle moves according to a random walk with identical mean and variance as the corresponding kinetic process. In the diffusive limit, where many collisions occur within a time step, the random walk corresponds to the diffusive limit of the Boltzmann-BGK equation. This diffusive step is meant to avoid the explicit simulation of a large number of collisions.

It is clear how this hybrid scheme solves the domain decomposition coupling issue: if, on average, less than one collision occurs in every time step, most of the particle trajectory will consist of kinetic steps, and the scheme corresponds to the kinetic approximation of the Boltzmann-BGK equation. If, on the other hand, more than one collision occurs in every time step, most of the particle trajectory will consist of diffusive steps, and the scheme corresponds to the diffusive approximation of the Boltzmann-BGK equation.

We will now briefly outline the details of the KDMC scheme from~\cite{mortier20}. Suppose a particle is released at time $t_\VAeventno$ with initial position $x(t_\VAeventno)=x_\VAeventno$ and velocity $v_k = \VAbackgroundspeedmean(x_\VAeventno)+\VAbackgroundspeedstdv(x_\VAeventno)\VAaddeventno{\VArnvel}{\VAeventno}$, where $\VArnvel_\VAeventno$ is a standard normal random number and $\VAbackgroundspeedmean(x_\VAeventno)$ and $\VAbackgroundspeedstdv^2(x_\VAeventno)$ are the mean and variance of the Maxwellian post-collisional velocity distribution at $x_\VAeventno$. The particle then moves with this constant velocity until a collision occurs. If the collision rate is given by $R(x)$, the time $\VAkintime_k$ until this collision is the solution of
\begin{equation}\label{eq:intinv}
\int_{0}^{\tau_k} R(x_k + v_k t) \mathrm{d}t = \epsilon_k
\end{equation}
where $\epsilon_k \sim \calE(1)$ is an exponentially distributed random number. Equivalently, using a change of variables $x=x_k+v_kt$, $\VAkintime_k$ is the solution of
\begin{equation}\label{eq:intinv-2}
\int_{x_k}^{x_k+v_k\tau_k} R(x) \mathrm{d}x = \epsilon_k v_k.
\end{equation}
In practice, and, also in our numerical experiments later on in~\secref{sec:numres}, the collision rate is such that $\tau_\VAeventno$ can easily be found from equation~\eqref{eq:intinv} or equation~\eqref{eq:intinv-2}, e.g., $R(x)$ is a piecewise constant or piecewise linear function.

The particle thus collides at time $t_\VAeventno+\tau_\VAeventno$, at a position $x_\VAeventno + v_\VAeventno\tau_\VAeventno$. In a standard kinetic simulation, the particle would now receive a new velocity from the Maxwellian post-collisional velocity distribution $v_{\VAeventno+1} = \VAbackgroundspeedmean(x_{\VAeventno+1})+\VAbackgroundspeedstdv(x_{\VAeventno+1})\VAaddeventno{\VArnvel}{{\VAeventno+1}}$, where $\VArnvel_{\VAeventno+1}$ is again a standard normal random number. The particle would then continue with this velocity until the next collision occurs. In the KD scheme, however, this new velocity $v_{k+1}$ is only applied after the particle moves diffusively for the remainder $\VAaddeventno{\VAdiftime}{\VAeventno}$ of the current time interval of length $\delta t$, where
\begin{equation}\label{eq:diffusive-time}
\VAaddeventno{\VAdiftime}{\VAeventno} = \delta t - (t_k + \tau_k \;\textrm{mod}\; \delta t).
\end{equation}
Thus, in a single KD step, the particle moves to a time
\begin{equation}\label{eq:time-update}
 t_{k+1} = t_{k} + \delta t = t_k+\tau_k+\theta_k,
\end{equation}
and a position
\begin{equation}\label{eq:position-update}
x_{k+1} = x_{k} + v_k\tau_k + a_k + d_k \chi_k,
\end{equation}
with $v_k\tau_k$ the kinetic contribution, and $a_k+d_k\chi_k$ the advection-diffusion contribution, and where $\chi_k\sim\calN(0,1)$ is a standard normal random number, and $a_k$ and $d_k$ represent the advective and diffusive part of the contribution, respectively. The advection-diffusion contribution is chosen such that the mean and variance of the random walk correspond exactly to the mean and variance of the actual kinetic process, conditioned such that the final velocity of the random walk is $\VAaddeventno{v}{{\VAeventno+1}}$, i.e., the velocity of the particle in the next time step, see~\cite{mortier20}. 

In the homogeneous case, when the quantities  $\VArate$, $\VAbackgroundspeedmean$, and $\VAbackgroundspeedstdv$ are constant, this can be achieved by choosing the advection coefficient as
\begin{equation}\label{eq:advection}
\VAaddeventno{a^\text{homo}}{\VAeventno}=\VAbackgroundspeedmean\VAaddeventno{\VAdiftime}{{\VAeventno}}+(v_{k+1}-\VAbackgroundspeedmean)\frac{1}{\VArate}\left(1-e^{-\VArate\VAaddeventno{\VAdiftime}{{\VAeventno}}}\right),
\end{equation}
and the diffusion coefficient as
\begin{equation}\label{eq:diffusion}
\VAaddeventno{d^\text{homo}}{\VAeventno}=\sqrt{\frac{2\VAbackgroundspeedstdv^2}{\VArate^2}\left(2(e^{-\VArate\VAaddeventno{\VAdiftime}{{\VAeventno}}}-1)+\VArate\VAaddeventno{\VAdiftime}{{\VAeventno}} (e^{-\VArate\VAaddeventno{\VAdiftime}{{\VAeventno}}}+1)\right)
+\frac{\left(v_{k+1}-\VAbackgroundspeedmean\right)^2}{\VArate^2}\left(1-2\VArate\VAaddeventno{\VAdiftime}{{\VAeventno}} e^{-\VArate\VAaddeventno{\VAdiftime}{{\VAeventno}}}-e^{-2\VArate\VAaddeventno{\VAdiftime}{{\VAeventno}}}\right)}\,,
\end{equation}
see~\cite{mortier20}. To cope with heterogeneity, i.e., when the collision rate $\VArate(x)$ and post-collisional velocity mean and variance $\VAbackgroundspeedmean(x)$ and $\VAbackgroundspeedstdv(x)$ are dependent on $x$, the advection and diffusion coefficient must be corrected to $\VAaddeventno{a^\text{hetero}}{\VAeventno}$ and $\VAaddeventno{d^\text{hetero}}{\VAeventno}$, respectively. The latter quantities differ from their homogeneous counterpart, because they use intermediate values for the mean and variance of the post-collision velocity evaluated at $x_{k+1/2}$, and because an additional advection term is added to cope with the case where $\VArate(x)$ varies strongly with $x$. The exact expressions for $\VAaddeventno{a^\text{hetero}}{\VAeventno}$ and $\VAaddeventno{d^\text{hetero}}{\VAeventno}$ are not shown here for the sake of breviety, but can be found entirely in~\cite{mortier20-b}. These adapted expressions are used in our numerical experiments in~\secref{sec:numres}.
 
The complete algorithm for a KD simulation is shown in~\algref{alg:KD}. The main procedure, \textsf{KineticDiffusion}, consists of a repeated call to the subroutine \textsf{KineticDiffusionStep}, where the latter implements a single KD step, as outlined above. The procedure returns $x_k$, the position of the particle at time $T$ after $k$ KD steps of length $\delta t$.

\begin{algorithm}
\begin{algorithmic}[1]
\Statex \textbf{input:} simulation end time $\VAendtime$, time step $\VAtimestep$
\Statex \textbf{output:} position $x_\VAeventno$ of the particle at the end time $\VAendtime$
\Statex
\Procedure{\textsf{KineticDiffusion}}{$\VAendtime, \VAtimestep$} \Comment{main KD routine}
\State $\VAeventno\gets0$
\State sample $\VArnvel_\VAeventno \sim \calN(0, 1)$ and $\VArnexp_\VAeventno \sim \calE(1)$
\State set $x_\VAeventno\gets1$, $v_\VAeventno\gets \VAbackgroundspeedmean + \VAbackgroundspeedstdv\VArnvel_\VAeventno$ and $t_\VAeventno\gets0$ \Comment{initialize position, velocity \& time}
\State solve equation~\eqref{eq:intinv} for $\VAkintime_\VAeventno$ using $\VArnexp_\VAeventno$ \Comment{compute time to first collision}
\While{$t_k + \tau_k < T$} \Comment{next collision is before end time $T$}
\State sample $\VArnvel_{\VAeventno+1} \sim \calN(0, 1)$, $\VArnexp_{\VAeventno+1} \sim \calE(1)$ and $\VArndif_\VAeventno \sim \calN(0, 1)$
\State $x_{\VAeventno+1}, v_{\VAeventno+1}, t_{\VAeventno+1}, \VAkintime_{\VAeventno+1}, \VAdiftime_{\VAeventno}\gets\mathsf{KineticDiffusionStep}(x_\VAeventno, v_\VAeventno, t_\VAeventno, \VAkintime_\VAeventno, \VArnvel_{\VAeventno+1}, \VArnexp_{\VAeventno+1}, \VArndif_\VAeventno)$
\State $\VAeventno \gets \VAeventno + 1$
\EndWhile
\If{$t_\VAeventno < T$} \Comment{move kinetically until end time $\VAendtime$}
\State $x_\VAeventno \gets x_\VAeventno +  v_\VAeventno(T - t_\VAeventno)$
\EndIf
\EndProcedure
\Statex
\Statex \textbf{input:} position $x_\VAeventno$, velocity $v_\VAeventno$ and current time $t_\VAeventno$ of the particle,
\Statex \hspace*{\WidthOfInput} time $\VAkintime_\VAeventno$ until the next kinetic event,
\Statex \hspace*{\WidthOfInput} standard normally distributed random number for the post-collisional velocity $\VArnvel_{\VAeventno+1}$, 
\Statex \hspace*{\WidthOfInput} exponentially distributed random number for the kinetic flight time $\VArnexp_{\VAeventno+1}$, and
 \Statex \hspace*{\WidthOfInput} standard normally distributed random number for the diffusive velocity $\VArndif_{\VAeventno}$
\Statex \textbf{output:} new position $x_{\VAeventno+1}$, new velocity $v_{\VAeventno+1}$ and new current time $t_{\VAeventno+1}$ of the particle,
\Statex \hspace*{\WidthOfOutput} time $\VAkintime_{\VAeventno+1}$ until the next kinetic event, duration $\VAdiftime_\VAeventno$ of previous diffusive step
\Statex
\Procedure{\textsf{KineticDiffusionStep}}{$x_\VAeventno, v_\VAeventno, t_\VAeventno, \VAkintime_{\VAeventno}, \VArnvel_{\VAeventno+1}, \VArnexp_{\VAeventno+1}, \VArndif_\VAeventno$} \Comment{KD subroutine}
\State $v_{\VAeventno+1} \gets \VAbackgroundspeedmean + \VAbackgroundspeedstdv\VArnvel_{\VAeventno+1}$ \Comment{new particle velocity in next kinetic phase}
\State $\VAdiftime_\VAeventno \gets \VAtimestep - (\VAkintime_\VAeventno\;\mathrm{mod}\;\VAtimestep)$ \Comment{duration of diffusive phase (equation~\eqref{eq:diffusive-time})}
\State compute $a_\VAeventno$ using~\eqref{eq:advection} and $d_\VAeventno$ using~\eqref{eq:diffusion}
\State $x_{\VAeventno+1} \gets x_\VAeventno + v_\VAeventno\VAkintime_\VAeventno + (a_\VAeventno + d_\VAeventno \VArndif_\VAeventno) $ \Comment{update of the particle position (equation~\eqref{eq:position-update})}
\State $t_{\VAeventno+1} \gets t_{\VAeventno} + \VAkintime_\VAeventno + \VAdiftime_\VAeventno$ \Comment{update of the particle time (equation~\eqref{eq:time-update})}
\State solve equation~\eqref{eq:intinv} for $\VAkintime_{\VAeventno+1}$ using $\VArnexp_{\VAeventno+1}$ \Comment{duration of next kinetic phase}
\EndProcedure
\end{algorithmic}
\caption{KD algorithm}
\label{alg:KD}
\end{algorithm}

The cost of~\algref{alg:KD} depends on the choice of the time step size $\delta t$. For large values of the time step size $\delta t$, the cost of the simulation depends on the number of time steps, since, on average and in a pure kinetic simulation, there will be many more collisions than time steps, but these many collisions are replaced by a single KD step. This cost continues to increase linearly with the number of time steps, until there is on average a single collision in every time interval. In that case, the cost of the KD scheme is equivalent to the cost of a purely kinetic simulation. For even smaller time step sizes, the simulation cost of~\algref{alg:KD} remains constant, since it amounts to recording each and every particle-background collision, and the average number of collisions is given by the collision rate, which is independent of the time step. Remark, however, that the result of a simulation with the KD scheme using a large time step is cheap but inaccurate: the simulation is biased. It is this bias that we will try to alleviate with the ML-KDMC scheme presented in~\secref{sec:mlapmc} below.

\section{Correlating fine and coarse particle paths\label{sec:ml_cor}}

A key component in the MLMC method is the ability to generate correlated samples for particle trajectories with different time step sizes. These correlated sample paths must be chosen such that they approximate the same underlying continuous particle trajectory, so that the difference between the simulated particle paths follows only from the difference in time step size, and not from the difference in kinetic behanviour. Below, we will present a new particle trajectory correlation scheme applied to the KD simulation scheme, which is one of the main contributions of this paper.

Let us introduce two time step sizes $\delta t_\level$ and $\delta t_{\level-1}$, with $\level$ the ``level'' of approximation. A particle trajectory with time step size $\delta t_\level$ thus corresponds to an approximation for the continuous particle trajectory with a relatively small time step size. This will be referred to as the \emph{fine} particle trajectory. Similarly, a particle trajectory with time step size $\delta t_{\level-1}$ corresponds to an approximation for the continuous particle trajectory with a relatively large time step size, and the trajectory will be referred to as the \emph{coarse} particle trajectory. We will add these level parameters $\level$ and $\level-1$ to all variables. For example, the initial position and velocity of the particle in the $k$th KD step with level $\level$ are denoted by $x_{\level, k}$ and $v_{\level, k}$, respectively. 

A good correlation between a fine sample particle trajectory at level $\level$ and a coarse sample particle trajectory at level $\level-1$ can be achieved by reusing the random numbers $\VArandparam_\level\coloneqq\{\nu_{\level,k}, \epsilon_{\level,k}, \chi_{\level,k}\}_{k=0}^{K_\level}$ from the simulation of the fine particle trajectory in the simulation of the coarse particle trajectory. To this end, let us define the operator $\varphi_{\level}^{\level-1} : \VArandparam_{\level} \mapsto \tilde\VArandparam_{\level_{j}}$, which maps the random numbers $\VArandparam_\level$ to a set of random numbers $\tilde\VArandparam_{\level-1}\coloneqq\{\tilde\nu_{\level-1,k}, \tilde\epsilon_{\level-1,k}, \tilde\chi_{\level-1,k}\}_{k=0}^{K_{\level-1}}$ for the simulation of the coarse particle trajectory. Applying the operator $\varphi_{\level}^{\level-1}$ consists of two phases. First, the random numbers on level $\level$ are mapped to the corresponding random numbers on level $\level-1$. Since the coarse particle trajectory corresponds to a simulation with a larger time step size, we expect fewer random numbers to be required at the coarse level, i.e., $K_\level > K_{\level-1}$, and hence more than one random number used for the fine particle trajectory must be mapped to a single random number from the coarse particle trajectory. We first need to determine which fine-level random numbers will be used for each coarse-level random number. We call this the \emph{mapping phase}. Afterwards, we must determine how the multiple random numbers from the fine particle trajectory should be aggregated into a single random number for the coarse particle trajectory. This is called the \emph{aggregation phase}. We will discuss these two phases in turn in~\secref{subsec:ml_cor_map}, respectively~\secref{subsec:ml_cor_agg}.

\begin{figure}
\centering
\begin{subfigure}{\textwidth}
\hspace{1.5cm}
	\tikzexternalenable
	\tikzsetnextfilename{random_number_mapping_standard}
\begin{tikzpicture}[xscale=3, every node/.style={outer sep=3pt}]

	\pgfmathsetmacro{\thetaZero}{.3} 
	\pgfmathsetmacro{\thetaOne}{.1} 
	\pgfmathsetmacro{\thetaTwo}{.2} 
	\pgfmathsetmacro{\thetaFour}{.3} 
	
	\foreach \i in {0, 1, ..., 4}{
		\draw[time step] (\i, -.1) -- (\i, .1);
		\node[anchor=south] at (\i, 0) {$\ifthenelse{\i=0}{0}{\ifthenelse{\i=1}{}{\i}{\VAtimestep_\level}}$};
	}
	\begin{scope}[on background layer]
		\kdstep{0}{\thetaZero}{1}
		\kdstep{1}{\thetaOne}{1}
		\kdstep{2}{\thetaTwo}{1}
		\kdstep{3}{\thetaFour}{1}
	\end{scope}
	\node[anchor=north] (k0f) at ($(0,0)!0.5!(\thetaZero,0)$) {$\VArnvel_{\level, 0}, \VArnexp_{\level, 0}$};
	\node[anchor=north] (d0f) at ($(\thetaZero,0)!0.5!(1,0)$) {$\VArndif_{\level, 0}$};
	\node[anchor=north] (k1f) at ($(1,0)!0.5!(1+\thetaOne,0)$) {$\VArnvel_{\level, 1}, \VArnexp_{\level, 1}$};
	\node[anchor=north] (d1f) at ($(1+\thetaOne,0)!0.5!(2,0)$) {$\VArndif_{\level, 1}$};
	\node[anchor=north] (k2f) at ($(2,0)!0.5!(2+\thetaTwo,0)$) {$\VArnvel_{\level, 2}, \VArnexp_{\level, 2}$};
	\node[anchor=north] (d2f) at ($(2+\thetaTwo,0)!0.5!(3,0)$) {$\VArndif_{\level, 2}$};
	\node[anchor=north] (k3f) at ($(3,0)!0.5!(3+\thetaFour,0)$) {$\VArnvel_{\level, 3}, \VArnexp_{\level, 3}$};
	\node[anchor=north] (d3f) at ($(3+\thetaFour,0)!0.5!(4,0)$) {$\VArndif_{\level, 3}$};
	
	\begin{scope}[yshift=-2cm]
		\foreach \i in {0, 1, ..., 2}{
			\draw[time step] (2*\i, -.1) -- (2*\i, .1);
			\node[anchor=north] at (2*\i, 0) {$\ifthenelse{\i=0}{0}{\ifthenelse{\i=1}{}{\i}{\VAtimestep_{\level-1}}}$};
		}
		\begin{scope}[on background layer]
			\kdstep{0}{\thetaZero}{2}
			\kdstep{2}{\thetaTwo}{2}
		\end{scope}
		\node[anchor=south] (k0c) at ($(0,0)!0.5!(\thetaZero,0)$) {$\tilde\VArnvel_{\level-1, 0}, \tilde\VArnexp_{\level-1, 0}$};
		\node[anchor=south] (d0c) at ($(\thetaZero,0)!0.5!(2,0)$) {$\tilde\VArndif_{\level-1, 0}$};
		\node[anchor=south] (k1c) at ($(2,0)!0.5!(2+\thetaTwo,0)$) {$\tilde\VArnvel_{\level-1, 1}, \tilde\VArnexp_{\level-1, 1}$};
		\node[anchor=south] (d1c) at ($(2+\thetaTwo,0)!0.5!(4,0)$) {$\tilde\VArndif_{\level-1, 1}$};
	\end{scope}
	
	\draw[default arrow] (k0f) -- (k0c);
	\draw[default arrow] (d0f) -- (d0c);
	\draw[default arrow] (k1f) -- (d0c);
	\draw[default arrow] (d1f) -- (d0c);
	\draw[default arrow] (k2f) -- (k1c);
	\draw[default arrow] (d2f) -- (d1c);
	\draw[default arrow] (k3f) -- (d1c);
	\draw[default arrow] (d3f) -- (d1c);
\end{tikzpicture}
	\tikzexternaldisable

\caption{The default case, where there is exactly one collision in every time interval of the fine particle trajectory. In this case, the index of the kinetic phase of the fine particle trajectory following the last diffusive phase that has already been aggregated is first $\kappa_0 = 0$, then $\kappa_1 = 2$, then $\kappa_2 = 4 \ldots$}
\label{fig:mapping-a}
\end{subfigure}
\begin{subfigure}{\textwidth}
\hspace{1.5cm}
	\tikzexternalenable
	\tikzsetnextfilename{random_number_mapping_timeshift}
\begin{tikzpicture}[xscale=3, every node/.style={outer sep=3pt}]

	\pgfmathsetmacro{\thetaZero}{.3} 
	\pgfmathsetmacro{\thetaOne}{1.5} 
	\pgfmathsetmacro{\thetaTwo}{.2} 
	
	\foreach \i in {0, 1, ..., 4}{
		\draw[time step] (\i, -.1) -- (\i, .1);
		\node[anchor=south] at (\i, 0) {$\ifthenelse{\i=0}{0}{\ifthenelse{\i=1}{}{\i}{\VAtimestep_\level}}$};
	}
	\begin{scope}[on background layer]
		\kdstep{0}{\thetaZero}{1}
		\kdstep{1}{\thetaOne}{2}
		\kdstep{3}{\thetaTwo}{1}
	\end{scope}
	\node[anchor=north] (k0f) at ($(0,0)!0.5!(\thetaZero,0)$) {$\VArnvel_{\level, 0}, \VArnexp_{\level, 0}$};
	\node[anchor=north] (d0f) at ($(\thetaZero,0)!0.5!(1,0)$) {$\VArndif_{\level, 0}$};
	\node[anchor=north] (k1f) at ($(1,0)!0.5!(1+\thetaOne,0)$) {$\VArnvel_{\level, 1}, \VArnexp_{\level, 1}$};
	\node[anchor=north] (d1f) at ($(1+\thetaOne,0)!0.5!(3,0)$) {$\VArndif_{\level, 1}$};
	\node[anchor=north] (k2f) at ($(3,0)!0.5!(3+\thetaTwo,0)$) {$\VArnvel_{\level, 2}, \VArnexp_{\level, 2}$};
	\node[anchor=north] (d2f) at ($(3+\thetaTwo,0)!0.5!(4,0)$) {$\VArndif_{\level, 2}$};
	
	\begin{scope}[yshift=-2cm]
		\foreach \i in {0, 1, ..., 2}{
			\draw[time step] (2*\i, -.1) -- (2*\i, .1);
			\node[anchor=north] at (2*\i, 0) {$\ifthenelse{\i=0}{0}{\ifthenelse{\i=1}{}{\i}{\VAtimestep_{\level-1}}}$};
		}
		\begin{scope}[on background layer]
			\kdstep{0}{\thetaZero}{2}
			\kdstep{2}{\thetaOne-1}{2}
		\end{scope}
		\node[anchor=south] (k0c) at ($(0,0)!0.5!(\thetaZero,0)$) {$\tilde\VArnvel_{\level-1, 0}, \tilde\VArnexp_{\level-1, 0}$};
		\node[anchor=south] (d0c) at ($(\thetaZero,0)!0.5!(2,0)$) {$\tilde\VArndif_{\level-1, 0}$};
		\node[anchor=south] (k1c) at ($(2,0)!0.5!(2+\thetaTwo,0)$) {$\tilde\VArnvel_{\level-1, 1}, \tilde\VArnexp_{\level-1, 1}$};
		\node[anchor=south] (d1c) at ($(2+\thetaTwo,0)!0.5!(4,0)$) {$\tilde\VArndif_{\level-1, 1}$};
	\end{scope}
	
	\draw[default arrow] (k0f) -- (k0c);
	\draw[default arrow] (d0f) -- (d0c);
	\draw[default arrow] (k1f) -- (k1c);
	\draw[default arrow] (d1f) -- (d1c);
	\draw[default arrow] (k2f) -- (d1c);
	\draw[default arrow] (d2f) -- (d1c);
\end{tikzpicture}
	\tikzexternaldisable

\caption{A case where no collision occurs during one time interval of the fine particle trajectory. In this case, the index of the kinetic phase of the fine particle trajectory following the last diffusive phase that has already been aggregated is first $\kappa_0 = 0$, then $\kappa_1 = 1$, then $\kappa_2 = 3 \ldots$}
\label{fig:mapping-b}
\end{subfigure}
\begin{subfigure}{\textwidth}
\hspace{1.5cm}
	\tikzexternalenable
	\tikzsetnextfilename{random_number_mapping_sum}
\begin{tikzpicture}[xscale=3, every node/.style={outer sep=3pt}]

	\pgfmathsetmacro{\thetaZero}{1.5} 
	\pgfmathsetmacro{\thetaTwo}{.3} 
	\pgfmathsetmacro{\thetaThree}{.1} 
	\pgfmathsetmacro{\thetaZeroCoarse}{2.3} 
	
	\foreach \i in {0, 1, ..., 4}{
		\draw[time step] (\i, -.1) -- (\i, .1);
		\node[anchor=south] at (\i, 0) {$\ifthenelse{\i=0}{0}{\ifthenelse{\i=1}{}{\i}{\VAtimestep_\level}}$};
	}
	\begin{scope}[on background layer]
		\kdstep{0}{\thetaZero}{2}
		\kdstep{2}{\thetaTwo}{1}
		\kdstep{3}{\thetaThree}{1}
	\end{scope}
	\node[anchor=north] (k0f) at ($(0,0)!0.5!(\thetaZero,0)$) {$\VArnvel_{\level, 0}, \VArnexp_{\level, 0}$};
	\node[anchor=north] (d0f) at ($(\thetaZero,0)!0.5!(2,0)$) {$\VArndif_{\level, 0}$};
	\node[anchor=north] (k1f) at ($(2,0)!0.5!(2+\thetaTwo,0)$) {$\VArnvel_{\level, 1}, \VArnexp_{\level, 1}$};
	\node[anchor=north] (d1f) at ($(2+\thetaTwo,0)!0.5!(3,0)$) {$\VArndif_{\level, 1}$};
	\node[anchor=north] (k2f) at ($(3,0)!0.5!(3+\thetaThree,0)$) {$\VArnvel_{\level, 2}, \VArnexp_{\level, 2}$};
	\node[anchor=north] (d2f) at ($(3+\thetaThree,0)!0.5!(4,0)$) {$\VArndif_{\level, 2}$};
	
	\begin{scope}[yshift=-2cm]
		\foreach \i in {0, 1, ..., 2}{
			\draw[time step] (2*\i, -.1) -- (2*\i, .1);
			\node[anchor=north] at (2*\i, 0) {$\ifthenelse{\i=0}{0}{\ifthenelse{\i=1}{}{\i}{\VAtimestep_{\level-1}}}$};
		}
		\begin{scope}[on background layer]
			\kdstep{0}{\thetaZeroCoarse}{4}
		\end{scope}
		\node[anchor=south] (k0c) at ($(0,0)!0.5!(2+\thetaTwo,0)$) {$\tilde\VArnvel_{\level-1, 0}, \tilde\VArnexp_{\level-1, 0}$};
		\node[anchor=south] (d0c) at ($(2+\thetaTwo,0)!0.5!(4,0)$) {$\tilde\VArndif_{\level-1, 0}$};
	\end{scope}
	
	\draw[default arrow] (k0f) -- (k0c);
	\draw[default arrow, shorten >=7mm] (d0f) -- (d0c);
	\draw[default arrow, shorten >=3mm] (k1f) -- (d0c);
	\draw[default arrow] (d1f) -- (d0c);
	\draw[default arrow] (k2f) -- (d0c);
	\draw[default arrow] (d2f) -- (d0c);
\end{tikzpicture}
	\tikzexternaldisable

\caption{A case where several kinetic and diffusive phases of the fine particle trajectory must be mapped onto a single diffusive phase of the coarse particle trajectory. In this case, the index of the kinetic phase of the fine particle trajectory following the last diffusive phase that has already been aggregated is first $\kappa_0 = 0$, then $\kappa_1 = 3, \ldots$}
\label{fig:mapping-c}
\end{subfigure}
\begin{subfigure}{\textwidth}
\hspace{1.5cm}
	\tikzexternalenable
	\tikzsetnextfilename{random_number_mapping_nefarious}
\begin{tikzpicture}[xscale=3, every node/.style={outer sep=3pt}]

	\pgfmathsetmacro{\thetaZero}{2.75} 
	\pgfmathsetmacro{\thetaOne}{1.2} 
	\pgfmathsetmacro{\thetaZeroCoarse}{.3} 
	
	\foreach \i in {0, 1, ..., 4}{
		\draw[time step] (\i, -.1) -- (\i, .1);
		\node[anchor=south] at (\i, 0) {$\ifthenelse{\i=0}{0}{\ifthenelse{\i=1}{}{\i}{\VAtimestep_\level}}$};
	}
	\begin{scope}[on background layer]
		\kdstep{0}{\thetaZero}{3}
		\kdstep{3}{\thetaOne}{\thetaOne}
	\end{scope}
	\node[anchor=north] (k0f) at ($(0,0)!0.5!(\thetaZero,0)$) {$\VArnvel_{\level, 0}, \VArnexp_{\level, 0}$};
	\node[anchor=north] (d0f) at ($(\thetaZero,0)!0.5!(3,0)$) {$\VArndif_{\level, 0}$};
	\node[anchor=north] (k1f) at ($(3,0)!0.5!(3+\thetaOne,0)$) {$\VArnvel_{\level, 1}, \VArnexp_{\level, 1}$};
	\node[anchor=north] (d1f) at ($(4.25,-.5)$) {};
	
	\begin{scope}[yshift=-2cm]
		\foreach \i in {0, 1, ..., 2}{
			\draw[time step] (2*\i, -.1) -- (2*\i, .1);
			\node[anchor=north] at (2*\i, 0) {$\ifthenelse{\i=0}{0}{\ifthenelse{\i=1}{}{\i}{\VAtimestep_{\level-1}}}$};
		}
		\begin{scope}[on background layer]
			\kdstep{0}{\thetaZeroCoarse}{2}
			\kdstep{2}{\thetaOne}{2}
		\end{scope}
		\node[anchor=south] (k0c) at ($(0,0)!0.5!(\thetaZeroCoarse,0)$) {$\tilde\VArnvel_{\level-1, 0}, \tilde\VArnexp_{\level-1, 0}$};
		\node[anchor=south] (d0c) at ($(\thetaZeroCoarse,0)!0.5!(2,0)$) {$\tilde\VArndif_{\level-1, 0}$};
		\node[anchor=south] (k1c) at ($(2,0)!0.5!(2+\thetaOne,0)$) {$\tilde\VArnvel_{\level-1, 1}, \tilde\VArnexp_{\level-1, 1}$};
		\node[anchor=south] (d1c) at ($(2+\thetaOne,0)!0.5!(4,0)$) {$\tilde\VArndif_{\level-1, 1}$};
	\end{scope}
	
	\draw[default arrow] (k0f) -- (k0c);
	\draw[default arrow, shorten >=4mm] (d0f) -- (d0c);
	\draw[default arrow] (k1f) -- (k1c);
	\draw[default arrow] (d1f) -- (d1c);
\end{tikzpicture}
	\tikzexternaldisable

\caption{A case of strong decorrelation, where no collision occurs during several time steps of the fine particle trajectory and both particles are in a part of the domain with a very different collision rate. In this case, the index of the kinetic phase of the fine particle trajectory following the last diffusive phase that has already been aggregated is first $\kappa_0 = 0$, then $\kappa_1 = 1$, then $\kappa_2 = 2 \ldots$}
\label{fig:mapping-d}
\end{subfigure}
\caption{Examples of the mapping of the random numbers from the fine to the coarse particle path.}
\label{fig:mapping}
\end{figure}

\subsection{Mapping the random numbers from the fine to the coarse particle path\label{subsec:ml_cor_map}}

We refer to~\figref{fig:mapping} for an illustration of the mapping phase. \figref{fig:mapping-a} shows the default case where there is exactly one collision in every KD time step of the fine particle trajectory. In this case, the kinetic phase of the coarse particle path uses the random numbers from the corresponding kinetic phase of the fine particle path, and the diffusive phase of the coarse particle path uses one diffusive, one kinetic and again one diffusive phase of the fine particle path. This mapping secures correlation by mapping the random numbers of the fine path to the coarse path based on the time they have an impact. \figref{fig:mapping-b} shows a case where there are time steps of the fine particle trajectory in which no collision occurs. In this case, we will only use the random numbers from the diffusive phase of the fine particle path for the diffusive phase of the coarse particle path. The length of the next kinetic phase of the coarse particle path, determined by $\epsilon_{\level, 2}$, will be adapted to reestablish the correlation, see the discussion of the aggregation phase in~\secref{subsec:ml_cor_agg}. \figref{fig:mapping-c} shows a case where several kinetic and diffusive phases of the fine particle trajectory must be mapped onto a single diffusive phase of the coarse particle trajectory. In this case, three diffusive phases and two kinetic phases must be aggregated to a single diffusive phase. Finally, \figref{fig:mapping-d} shows how the mapping would behave in a quite dramatic situation in which the collision rate of the coarse particle is much larger than the collision rate of the fine particle. Such a situation can occur only if the fine and coarse particle trajectories are strongly decorrelated and the coarse and fine particles are in different parts of the domain with very different collision rates. We did not observe such a strong decorrelation in our numerical results in~\secref{sec:numres}.

The mapping can be described mathematically as follows. Suppose that $k-1$ KD steps have already been taken in the coarse particle path. During these steps, several kinetic and diffusive parts of the fine particle trajectory have already been used in the aggregation. Let $\VAfineindex{\VAeventno}$ be the index of the first kinetic phase of the fine particle trajectory following the last diffusive phase that has already been aggregated, see~\figref{fig:mapping}. The kinetic phase of the $k$th KD step in the coarse particle trajectory then uses the random numbers $\nu_{\level, \kappa_k}$ and $\epsilon_{\level, \kappa_k}$, i.e., the random numbers of the $\kappa_k$th kinetic phase of the fine particle trajectory. Recall that $\nu$ and $\epsilon$ determine the velocity and length of the kinetic phase, respectively. The diffusive phase of the $k$th KD step in the coarse particle trajectory combines the motion of all the diffusive and kinetic phases of the fine particle trajectory up to but not including the next kinetic phase of which the collision takes place after the end of the current coarse time step $\delta t_{\level-1}$. The index of this last kinetic phase is then $\kappa_{k+1}$. That way, the random numbers of the fine trajectory are mapped to the coarse trajectory in such a way that the time during which a random variable has an effect on the particle trajectory is approximately equal for the fine and the coarse particle. The random numbers used in the $k$th coarse diffusive phase are thus an aggregation of the random numbers $\nu_{\level, m}$ and $\epsilon_{\level, m}$,$m=\kappa_k+1, \ldots, \kappa_{k+1}-1$, and $\chi_{\level, m}$, $m=\kappa_k, \ldots, \kappa_{k+1}-1$, used in the fine particle trajctory.

\subsection{Aggregating the random numbers in the coarse path\label{subsec:ml_cor_agg}}

After the mapping phase, outlined in~\secref{subsec:ml_cor_map}, we must now specify \emph{how} the random numbers $\VArandparam_\level$ from the fine particle trajectory will be aggregated to the random numbers $\tilde\VArandparam_{\level-1}$ for the coarse particle trajectory. The aggregation of the kinetic phase is straightforward. In the $k$th kinetic phase of the coarse particle trajectory, the random variable for the new Maxwellian post-collisional velocity is simply set to
\begin{equation}\label{eq:nu-update}
\tilde\nu_{\level-1, k} \coloneqq \nu_{\level, \kappa_k},
\end{equation}
where we recall the notation $\kappa_k$ for the index of the next unused kinetic phase of the fine particle trajectory after $k-1$ KD steps in the coarse particle path. The length of the kinetic path is determined by
\begin{equation}\label{eq:epsilon-update}
\tilde\epsilon_{\level-1, k} \coloneqq \epsilon_{\level, \kappa_k} - \int_0^{t_{\level-1, k} - t_{\level, \kappa_k}} R(x_{\level, \kappa_k} + v_{\level, \kappa_k} t) \mathrm{d}t,
\end{equation}
see~\figref{fig:mapping-b}. Because of the memorylessness of the exponential distribution, the random variable $\tilde\epsilon_{\level-1, k}$ still follows the required (exponential) distribution.

In the $k$th diffusive phase of the coarse particle trajectory, the random number aggregation for $\tilde\chi_{\level-1, k}$ will use a weighted sum of the normally distributed numbers $\VArnvel_{\level, m}$ with weights $\alpha_m$, not all zero, for $m=\kappa_k+1, \ldots, \kappa_{k+1}-1$, and normally distributed numbers $\VArndif_{\level, m}$ with weights $\beta_m$, not all zero, for $m=\kappa_k, \ldots, \kappa_{k+1}-1$. With a suitable normalization, we set
\begin{equation}\label{eq:chi-update}
\tilde\chi_{\level-1, k} \coloneqq \frac{\beta_{\kappa_k}\chi_{\level, \kappa_k} + \displaystyle\sum_{m=\kappa_k+1}^{\kappa_{k+1}-1}\left(\alpha_m \nu_{\level, m} + \beta_m \chi_{\level, m}\right)}{\sqrt{\beta_{\kappa_k}^2 + \displaystyle\sum_{m=\kappa_k+1}^{\kappa_{k+1}-1} \left(\alpha_m^2 + \beta_m^2 \right)}}.
\end{equation}
See~\figref{fig:mapping} for an illustration of the aggregation.

To ensure a good correlation between the fine and coarse particle trajectory, the weights $\alpha_m$ and $\beta_m$ are chosen such that the aggregated random numbers for the kinetic velocity $\nu_{\level, m}$ and the diffusive velocity $\chi_{\level, m}$ of the fine particle trajectory have an approximately equal effect on the diffusive velocity $\chi_{\level-1, m}$ of the coarse particle trajectory. We recall that the position of a particle following the fine trajectory can be computed as
\begin{equation}
x_{\level, K_\level} = \underbrace{\VAaddeventno{\VAkintime}{{\level,0}}(\VAbackgroundspeedmean+\VAbackgroundspeedstdv\VAaddeventno{\VArnvel}{{\level, 0}})}_\text{kinetic}+\underbrace{\VAaddeventno{a}{{\level,0}}+\VAaddeventno{d}{{\level,0}}\VAaddeventno{\VArndif}{{\level,0}}}_\text{diffusion}+\underbrace{\VAaddeventno{\VAkintime}{{\level,1}}(\VAbackgroundspeedmean+\VAbackgroundspeedstdv\VAaddeventno{\VArnvel}{{\level,1}})}_\text{kinetic}+\underbrace{\VAaddeventno{a}{{\level,1}}+\VAaddeventno{d}{{\level,1}}\VAaddeventno{\VArndif}{{\level,1}}}_\text{diffusion}+\cdots\,,\label{eq:ml_cor_finesum}
\end{equation}
where $\VAaddeventno{\VAkintime}{{\level, k}}$ is determined by $\epsilon_{\level, k}$, see equation~\eqref{eq:intinv}, and where $a_{\level, k}$ and $d_{\level, k}$ depend on $\epsilon_{\level, k}$ and $\nu_{\level, k+1}$, see equations~\eqref{eq:advection} and~\eqref{eq:diffusion}. In effect, this means that the random variable $\nu_{\level, k+1}$ for the Maxwellian post-collisional velocity plays a role for a longer time than just the duration of the kinetic phase $\tau_{\level, k+1}$, since it appears in the expression for the advective and diffusive contributions $a_{\level, k}$, respectively, $d_{\level, k}$, of the previous KD phase. This excess time thus also appears in the duration of the previous diffusive phase $\theta_{\level, k}$. To incorporate this aspect of dependence between different time steps in our aggregation scheme, we transfer the expected time during which the random number $\nu_{\level, k+1}$ for the post-collisional velocity plays a role in the diffusive phase to the next kinetic phase. This time, which we will denote by $\kiessymbool_{\level,m}$, is equal to the expected time of the diffusive phase in which this final velocity is used, i.e.,
\begin{equation}
\kiessymbool_{\ell,m}\coloneqq\frac{1}{R_{m}}\left(1-e^{-R_{m}\theta_{\level,{m}}}\right).
\end{equation}
The duration of the kinetic and diffusive phase are thus updated to
\begin{equation}
\tau_{\ell,m+1}'\coloneqq\tau_{\ell,m+1}+\kiessymbool_{\ell,m}\, , \quad \text{respectively} \quad \theta_{\ell,m}' \coloneqq \theta_{\ell,m}-\kiessymbool_{\ell,m}.
\end{equation}

The weights $\alpha_m$, $m=\kappa_k+1, \ldots, \kappa_{k+1}-1$ for the normally distributed random numbers $\nu_{\level, m}$ are chosen as the modified kinetic time during which the velocity $\nu_{\ell,m}$ determines the kinetic motion, i.e.,
\begin{equation}
\VAaddeventno{\alpha}{m}\coloneqq\VAaddeventno{\VAkintime'}{{\level,m}}\VAbackgroundspeedstdv.
\end{equation}

The weights $\beta_m$, $m=\kappa_k, \ldots, \kappa_{k+1}-1$ for the normally distributed random numbers $\chi_{\level, m}$ are chosen as
\begin{equation}\label{eq:weights-beta}
\VAaddeventno{\beta}{m}\coloneqq\sqrt{\frac{2\VAbackgroundspeedstdv^2}{\VArate_m^2}\left(e^{-\VArate_m\VAaddeventno{\VAdiftime'}{{\level,m}}}+\VArate_m\VAdiftime'_{\level, m} - 1\right)}.
\end{equation}
These weights are equal to the standard deviation of the diffusive position update without conditioning on the final velocity $v_{\level, {m+1}}$.

\subsection{An algorithm for correlated KD sampling}

A complete algorithm for correlated KD sampling is shown in~\algref{alg:CorrelatedKD}. The routine \textsf{CorrelatedKineticDiffusion} moves the particle simultaneously on the coarse and fine particle trajectory in one sweep over the time domain. The procedure reuses the subroutine \textsf{KineticDiffusionStep} from~\algref{alg:KD}. The output of the procedure is $x_{\level, k_1}$ and $x_{\level-1, k_2}$, two correlated samples of the particle positions at time $T$, approximated with a time step $\delta t_\level$, respectively $\delta t_{\level-1}$.

\begin{algorithm}
\begin{algorithmic}[1]
\Statex \textbf{input:} simulation end time $\VAendtime$, time steps $\VAtimestep_{\level}$ and $\VAtimestep_{\level-1}$
\Statex \textbf{output:} positions $x_{\level, k_1}$ and $x_{\level-1, k_2}$ of the particle at the end time $\VAendtime$
\Statex
\Procedure{\textsf{CorrelatedKineticDiffusion}}{$\VAendtime, \VAtimestep_\level, \VAtimestep_{\level-1}$}
\State set $k_1\gets0$ and $k_2\gets0$
\State sample $\VArnvel_{\level, k_1} \sim \calN(0, 1)$ and $\VArnexp_{\level, k_1} \sim \calE(1)$
\State set $\tilde\VArnvel_{\level-1, k_2}\gets\VArnvel_{\level, k_1}$ and $\tilde\VArnexp_{\level-1, k_2}\gets\VArnexp_{\level, k_1}$
\State set $x_{\level, k_1}\gets1$, $v_{\level, k_1}\gets \VAbackgroundspeedmean + \VAbackgroundspeedstdv\VArnvel_{\level, k_1}$ and $t_{\level,k_1}\gets0$ \Comment{initialize position, velocity \& time}
\State set $x_{\level-1, k_2}\gets1$, $v_{\level-1, k_2}\gets \VAbackgroundspeedmean + \VAbackgroundspeedstdv\tilde\VArnvel_{\level-1,k_2}$ and $t_{\level-1,k_2}\gets0$ \Comment{... same for coarse path}
\State solve equation~\eqref{eq:intinv} for $\VAkintime_{\level, k_1}$ and set $\VAkintime_{\level-1, k_2}\gets\VAkintime_{\level, k_1}$ \Comment{compute time to first collision}
\While{$t_{\level-1} + \tau_{\level-1} < T$ or $t_{\level} + \tau_{\level} < T$} \Comment{next coarse collision is before end time $T$}
\LineComment{repeat until next random number aggregation from fine to coarse path}
\Repeat \Comment{move the fine particle path}
\State sample $\VArnvel_{\level, k_1+1} \sim \calN(0, 1)$, $\VArnexp_{\level, k_1+1} \sim \calE(1)$ and $\VArndif_{\level, k_1} \sim \calN(0, 1)$
\State $x_{\level, k_1+1}, v_{\level, k_1+1}, t_{\level, k_1+1}, \VAkintime_{\level, k_1+1}, \VAdiftime_{\level,k_1}\gets$
\State \hfill $\mathsf{KineticDiffusionStep}(x_{\level,k_1}, v_{\level,k_1}, t_{\level,k_1}, \VAkintime_{\level,k_1}, \VArnvel_{\level,k_1+1}, \VArnexp_{\level,k_1+1}, \VArndif_{\level,k_1})$
\State $k_1 \gets k_1 + 1$
\Until{$\lceil (t_{\level,k_1} + \VAkintime_{\level,k_1}) / \VAtimestep_\level \rceil \VAtimestep_\level > \lceil (t_{\level-1,k_2} + \VAkintime_{\level-1,k_2}) / \VAtimestep_{\level-1} \rceil \VAtimestep_{\level-1}$}
\If{$t_{\level-1} + \tau_{\level-1} < T$}
\State set $\tilde\VArnvel_{\level-1, k_2+1} \gets \VArnvel_{\level, k_1}$ \Comment{equation~\eqref{eq:nu-update}}
\State compute $\tilde\VArnexp_{\level-1, k_2+1}$ according to equation~\eqref{eq:epsilon-update} and $\tilde\VArndif_{\level-1, k_2}$ according to equation~\eqref{eq:chi-update}
\State $x_{\level-1, k_2+1}, v_{\level-1, k_2+1}, t_{\level-1, k_2+1}, \VAkintime_{\level-1, k_2+1},  \VAdiftime_{\level-1,m}\gets$
\State \hfill$\mathsf{KineticDiffusionStep}(x_{\level-1,k_2}, v_{\level-1,k_2}, t_{\level-1,k_2}, \VAkintime_{\level-1,k_2}, \tilde\VArnvel_{\level-1,k_2+1}, \tilde\VArnexp_{\level-1,k_2+1}, \tilde\VArndif_{\level-1,k_2})$
\State $k_2 \gets k_2 + 1$
\EndIf
\EndWhile
\State $k_1\gets\textsf{min}(k_1|t_{\level,k_1}+\tau_{\level,k_1}\geq \VAendtime)$\Comment{remove fine steps beyond $\VAendtime$}
\If{$t_{\level,k_1} < T$} \Comment{move fine particle path kinetically until time $\VAendtime$}
\State $x_{\level,k_1} \gets x_{\level,k_1} +  v_{\level,k_1}(T - t_{\level,k_1})$
\EndIf
\If{$t_{\level-1,k_2} < T$} \Comment{move coarse particle path kinetically until time $\VAendtime$}
\State $x_{\level-1,k_2} \gets x_{\level-1,k_2} +  v_{\level-1,k_2}(T-t_{\level-1,k_2})$
\EndIf
\EndProcedure
\end{algorithmic}
\caption{Correlated KD algorithm}
\label{alg:CorrelatedKD}
\end{algorithm}

\section{Multilevel Kinetic-Diffusion Monte Carlo}\label{sec:mlapmc}

In this section, we outline the details of our ML-KDMC method. First, in~\secref{subsec:mlmc}, we briefly discuss the MLMC method. This multilevel method uses a hierarchy of coarser approximations with ever larger time step sizes, to reduce the cost of the MC simulation. Next, in~\secref{subsec:levelselect}, we discuss the specific challenges in applying the MLMC method to the KDMC scheme from~\secref{sec:ml_sim}. Notably, the behanviour of variance and cost of the multilevel differences with increasing level parameter $\level$ is not the usual monotone behanviour as observed in models described by differential equations, see~\cite{cliffe11, giles08}. This atypical behanviour poses challenges for the level selection problem, i.e., the choice of the coarser time step sizes that are included in the multilevel hierarchy.

\subsection{Multilevel sampling}\label{subsec:mlmc}

Let us introduce the set of time step sizes $\{\VAtimestep_\level \coloneqq \VAendtime / 2^\level\}_{\level=0}^L$. These time step sizes discretize the time domain $[0, \VAendtime]$ into $2^\level$ time intervals $[\VAtimestepno\cdot\VAtimestep_\level, (\VAtimestepno+1)\cdot\VAtimestep_\level]$ for $\VAtimestepno=0, 1, \ldots, 2^\level-1$. Furthermore, let $x_\level(t, \VArandparam_\level)$ denote the position of a particle at time $t$, computed using a time step $\VAtimestep_\level=\VAendtime/2^\level$, where $\VAendtime$ is the end time and $\level\ge0$ is the level of approximation. Here, we explicitly denote the dependence of the particle position $x_\level$ on the random variables $\VArandparam_\level = \{\VArnvel_{\level,\VAeventno}, \VArnexp_{\level,\VAeventno}, \VArndif_{\level,\VAeventno}\}_{\VAeventno=0}^{\VAncollisions_\level}$ used to simulate the trajectory of the particle. Our goal is to compute the position of the particle at the end time $\VAendtime$, i.e., $x_\level(\VAendtime, \VArandparam_\level)$. This position is the \emph{quantity of interest}, hereafter denoted as
\begin{equation*}
\VAqoi_\level(\VArandparam_\level) \coloneqq x_\level(\VAendtime, \VArandparam_\level).
\end{equation*}
Sometimes, the argument will be dropped when the meaning is clear from the context, i.e., we write $\VAqoi_\level$ instead of $\VAqoi_\level(\VArandparam_\level)$. The $n$th independent identically distributed sample of this random variable $\VAqoi_\level$ is denoted as
\begin{equation*}
\VAqoi_\level(\VArandparam_\level^{(\VAparticleno)}) \coloneqq x_\level(\VAendtime, \VArandparam_\level^{(\VAparticleno)}) \text{ with } \VArandparam_\level^{(\VAparticleno)} \coloneqq \{\VArnvel_{\level, \VAeventno}^{(\VAparticleno)}, \VArnexp_{\level, \VAeventno}^{(\VAparticleno)}, \VArndif_{\level, \VAeventno}^{(\VAparticleno)}\}_{k=0}^{\VAncollisions_\level^{(\VAparticleno)}}.
\end{equation*}
Furthermore, let $\{0,1,\ldots,L\}$ be the set of all possible levels, in increasing order, and let $\calL \subseteq \{0, 1, \ldots, L\}$ be an ordered subset of $J$ levels, with $0<J\leq L+1$, defined as $\calL\coloneqq\{\level_j\}_{j=1}^J$ where $0 \leq   \level_1 < \ldots < \level_J  = L$. Note that we restrict our attention to subsets that include the level with highest accuracy $L$, i.e., $\level_J = L$. If $J=1$, then $\calL=\{L\}$, and only the most accurate level $\level_J=L$ is used. On the other hand, if $J=L+1$, then all levels $0, 1, 2, \ldots, L$ are included in the set $\calL$. The set of all feasible subsets $\calL$ will be denoted by $S_\calL$.

The MLMC estimator for the expected value of the quantity of interest on level $L$ is
\begin{equation}\label{eq:mlmc}
\VAqoiest_{\calL} \coloneqq \frac{1}{N_{\level_1}} \sum_{n=1}^{N_{\level_1}} \VAqoi_{\level_1}(\VArandparam_{\level_1}^{(\VAparticleno)}) + \sum_{j=2}^J \frac{1}{N_{\level_j}} \sum_{n=1}^{N_{\level_j}} \left(\VAqoi_{\level_j}(\VArandparam_{\level_j}^{(\VAparticleno)}) - \VAqoi_{\level_{j-1}}(\tilde\VArandparam_{\level_{j-1}}^{(\VAparticleno)}) \right),
\end{equation}
where all $\VArandparam_{\level_j}^{(\VAparticleno)}$ are independent and identically distributed samples for each $n=1, 2, \ldots, N_{\level_j}$ and for each $\level_j \in \calL$, with $N_{\level_j}$ the total number of samples on level $\level_j$. The random variables $\tilde\VArandparam_{\level_{j-1}}^{(\VAparticleno)} $ are defined as
\begin{equation}
\tilde\VArandparam_{\level_{j-1}}^{(\VAparticleno)} \coloneqq \{\tilde\VArnvel_{\level_{j-1}, \VAeventno}^{(\VAparticleno)}, \tilde\VArnexp_{\level_{j-1}, \VAeventno}^{(\VAparticleno)}, \tilde\VArndif_{\level_{j-1}, \VAeventno}^{(\VAparticleno)}\}_{\VAeventno=0}^{\tilde\VAncollisions_{\level_{j-1}}^{(\VAparticleno)}},
\end{equation}
where $\tilde\VArnvel_{\level_{j-1}, \VAeventno}^{(\VAparticleno)}$, $\tilde\VArnexp_{\level_{j-1}, \VAeventno}^{(\VAparticleno)}$ and $\tilde\VArndif_{\level_{j-1}, \VAeventno}^{(\VAparticleno)}$ are computed using the mapping and aggregation $\varphi_{\level_j}^{\level_{j-1}} : \VArandparam_{\level_j}^{(\VAparticleno)} \mapsto \tilde\VArandparam_{\level_{j-1}}^{(\VAparticleno)}$ from~\secref{subsec:ml_cor_agg} applied to $\VArandparam_{\level_{j}}^{(\VAparticleno)}$, for each $\VAeventno=0, 1, \ldots, \tilde\VAncollisions_{\level_{j-1}}^{(\VAparticleno)}$, and with $\tilde\VAncollisions_{\level_{j-1}}^{(\VAparticleno)}$ the number of collisions in the $n$th coarse aggregated particle path. 

Using the shorthand notation $\Delta \VAqoi_{\level_i, \level_j} \coloneqq \VAqoi_{\level_i}(\VArandparam_{\level_i}) - \VAqoi_{\level_j}(\tilde\VArandparam_{\level_j})$ for a multilevel \emph{difference}, the MLMC estimator can be written compactly as
\begin{equation}\label{eq:mlmc-compact}
\VAqoiest_\calL = \frac{1}{N_{\level_1}}  \sum_{\VAparticleno=1}^{N_{\level_1}} \VAqoi_{\level_1}^{(n)} + \sum_{j=2}^J \frac{1}{N_{\level_j}} \sum_{\VAparticleno=1}^{N_{\level_j}} \Delta \VAqoi_{\level_j, \level_{j-1}}^{(\VAparticleno)},
\end{equation}
where $\Delta \VAqoi_{\level_j, \level_{j-1}}^{(\VAparticleno)}\coloneqq\VAqoi_{\level_j}(\VArandparam_{\level_j}^{(\VAparticleno)}) - \VAqoi_{\level_{j-1}}(\tilde\VArandparam_{\level_{j-1}}^{(\VAparticleno)})$ denotes the $\VAparticleno$th realization of $\Delta \VAqoi_{\level_j, \level_{j-1}}$, and $\VAqoi_{\level_1}^{(n)} \coloneqq \VAqoi_{\level_1}(\VArandparam_{\level_1}^{(\VAparticleno)})$. For later use, let us introduce the following notation. We define $E_{\level_i} \coloneqq |\E{\VAqoi_{\level_i}}|$ and $E_{\level_i,\level_j} \coloneqq |\E{\Delta \VAqoi_{\level_i,\level_j}}|$, where $\E{\;\cdot\;}$ denotes the expected value operator. Similarly, we write $V_{\level_i} \coloneqq \V{\VAqoi_{\level_i}}$ and $V_{\level_i,\level_j} \coloneqq \V{\Delta \VAqoi_{\level_i,\level_j}}$, where $\V{\;\cdot\;}$ denotes the variance operator. Finally, we define  $C_{\level_i} \coloneqq \cost{\VAqoi_{\level_i}}$ and $C_{\level_i,\level_j} \coloneqq \cost{\VAqoi_{\level_i}} + \cost{\VAqoi_{\level_j}}$, where $\cost{\;\cdot\;}$ denotes the cost of computing a single realization of a random variable. This cost can be measured in, for example, the number of floating point operations, or in actual wall clock time. We will use the latter in our numerical experiments.

It is easy to see that the multilevel estimator is an unbiased estimator for $\E{\VAqoi_L}$, since
\begin{align}
\E{\VAqoiest_{\calL}} &= \E{ \VAqoi_{\level_1}(\VArandparam_{\level_1})} + \sum_{j=2}^J \E{\VAqoi_{\level_j}(\VArandparam_{\level_j}) - \VAqoi_{\level_{j-1}}(\tilde\VArandparam_{\level_{j-1}})} \\
&= \vphantom{\sum_{j=2}^J}\E{ \VAqoi_{\level_1}(\VArandparam_{\level_1})} + \E{\VAqoi_{\level_2}(\VArandparam_{\level_2}) - \VAqoi_{\level_1}(\tilde\VArandparam_{\level_1})}
+ \ldots + \E{\VAqoi_{\level_J}(\VArandparam_{\level_J}) - \VAqoi_{\level_{J-1}}(\tilde\VArandparam_{\level_{J-1}})} \label{eq:telsum}\\
&= \vphantom{\sum_{j=2}^J}\E{\VAqoi_{\level_J}(\VArandparam_{\level_J}) } = \E{\VAqoi_{L} }, \label{eq:unbiased}
\end{align}
where we used the linearity of the expectation operator and the fact that the random variables $\VArandparam_{\level_{j}}$ and $\tilde\VArandparam_{\level_{j-1}}$ have the same distribution. We will numerically verify this assumption later on in our experiments in~\secref{subsec:consistency}. Note that equation~\eqref{eq:unbiased} is the motivation for our earlier restriction that $\ell_J=L$. By including the level $L$ with highest accuracy into the set $\calL$, the multilevel estimator is unbiased estimator for the expected value of the quantity of interest on that level. The variance of the multilevel estimator can be expressed as
\begin{equation}
\V{\VAqoiest_{\calL}} = \frac{\V{ \VAqoi_{\level_1}(\VArandparam_{\level_1})}}{N_{\level_1}} + \sum_{j=2}^J \frac{\V{\VAqoi_{\level_j}(\VArandparam_{\level_j}) - \VAqoi_{\level_{j-1}}(\tilde\VArandparam_{\level_{j-1}})}}{N_{\level_j}}
= \frac{V_{\level_1}}{N_{\level_1}} +\sum_{j=2}^J \frac{V_{\level_j, \level_{j-1}}}{N_{\level_j}},
\end{equation}
where we again used the observation that the random variables $\VArandparam_{\level_{j}}$ and $\tilde\VArandparam_{\level_{j-1}}$ have the same distribution, and that $\bsy_{\level_j}$ and $\bsy_{\level_k}$ are independent for any $j \ne k$.

There are two sources of error in the MLMC estimator in equation~\eqref{eq:mlmc-compact}: the discretization error, related to the finite time step $\VAtimestep_L$, and the statistical error, present because we replace the expected value by a sample average of a finite set of samples. The accuracy of the estimator can be quantified using the mean square error (MSE), where these two sources of error become apparent:
\begin{align}
\mse{\VAqoiest_{\calL}} &\coloneqq \E{\left(\VAqoiest_{\calL} - \E{\VAqoi}\right)^2} \\
&= \left(\E{\VAqoiest_{\calL}} - \E{\VAqoi}\right)^2 + \E{\left(\VAqoiest_{\calL} - \E{\VAqoiest_{\calL}}\right)^2} \\
&=  (\underbrace{\vphantom{\sum_{j=2}^J}\E{\VAqoi_L-\VAqoi}}_{\textstyle\strut \text{bias}})^2 + \Bigg( \underbrace{\frac{V_{\level_1}}{N_{\level_1}} +\sum_{j=2}^J \frac{V_{\level_j, \level_{j-1}}}{N_{\level_j}}}_{\textstyle\strut \text{variance}}\Bigg)\label{eq:mse}
\end{align}
The first term in equation~\eqref{eq:mse} is the square of the bias of the estimator, representing the discretization error. The second term in equation~\eqref{eq:mse} is the variance of the estimator, representing the stochastic part of the error. To impose an MSE of at most $\varepsilon^2$, or, equivalently, a \emph{root mean square error} (RMSE) of at most $\varepsilon$, it is now sufficient to enforce that
\begin{align}
|\E{\VAqoi_L-\VAqoi}| &\le \frac{\varepsilon}{\sqrt{2}} &&\text{(bias constraint), and} \label{eq:bias_constraint} \\
\frac{V_{\level_1}}{N_{\level_1}} +\sum_{j=2}^J \frac{V_{\level_j, \level_{j-1}}}{N_{\level_j}} &\le \frac{\varepsilon^2}{2} &&\text{(statistical constraint)}.\label{eq:statistical_contraint}
\end{align}

Two unknowns remain in the formulation of the MLMC estimator in equation~\eqref{eq:mlmc-compact}.
\begin{itemize}
\item The \textbf{choice of levels} $\calL = \{\ell_j\}_{j=1}^J$. The optimal set of levels can be found by minimizing the total cost of the estimator, while ensuring that the bias constraint, i.e., equation~\eqref{eq:bias_constraint}, is satisfied.
\item The \textbf{number of samples} $N_{\level_j}$, for each $j=1, 2, \ldots, J$. Given the set of levels $\calL$, this number of samples can be found by minimizing the total cost of the estimator, while ensuring that the statistical constraint, i.e., equation~\eqref{eq:statistical_contraint}, is satisfied.
\end{itemize}
We will address the latter problem, i.e., determining the optimal number of samples $\{N_{\level_j}\}_{j=1}^J$, in the remainder of this section. The first problem, i.e., the level selection strategy, will be tackled in the next section. The reason for this turnaround is that easy-to-obtain analytic expressions exist for the optimal number of samples $N_{\level_j}$ at each level $\level_j$.

The total cost of the MLMC estimator can be expressed as
\begin{equation}
\cost{\VAqoiest_{\calL}} = N_{\level_1} C_{\level_1} + \sum_{j=2}^J N_{\level_j} C_{\level_j, \level_{j-1}},
\end{equation}
where we recall the notation $C_{\level_1}$ for the cost of computing a single realization of $Q_{\level_1}$, and $C_{\level_j, \level_{j-1}}$ for the cost of computing a single realization of the difference $\Delta Q_{\level_j, \level_{j-1}}$. Next, consider the constrained minimization problem
\begin{align}
\min_{N_{\level_1}, N_{\level_2}, \ldots, N_{\level_J}} &N_{\level_1} C_{\level_1} + \sum_{j=2}^J N_{\level_j} C_{\level_j, \level_{j-1}} \\
\text{s.t.}\qquad &\frac{V_{\level_1}}{N_{\level_1}} +\sum_{j=2}^J \frac{V_{\level_j, \level_{j-1}}}{N_{\level_j}} = \frac{\varepsilon^2}{2},
\end{align}
see, e.g.,~\cite{giles08} for details. The Lagrangian of this problem is
\begin{equation}
\mathscr{L}(N_{\level_1}, N_{\level_2}, \ldots, N_{\level_J}) = N_{\level_1} C_{\level_1} + \sum_{j=2}^J N_{\level_j} C_{\level_j, \level_{j-1}} + \zeta \left(\frac{V_{\level_1}}{N_{\level_1}} +\sum_{j=2}^J \frac{V_{\level_j, \level_{j-1}}}{N_{\level_j}} - \frac{\varepsilon^2}{2}\right),
\end{equation}
where $\zeta$ is a Lagrange multiplier and where we treat the unknowns $N_{\level_1}$, $N_{\level_2}$, $\ldots$, $N_{\level_J}$ as continuous variables. Proceeding as usual, the first-order necessary optimality conditions are
\begin{equation}
\begin{cases}
\begin{aligned}
\displaystyle\frac{\partial \mathscr{L}}{\partial N_{\level_1}} &=  C_{\level_1} - \zeta \displaystyle\frac{V_{\level_1}}{N_{\level_1}^2} = 0,\\
\displaystyle\frac{\partial \mathscr{L}}{\partial N_{\level_j}} &=  C_{\level_j, \level_{j-1}} - \zeta \displaystyle\frac{V_{\level_j, \level_{j-1}}}{N_{\level_j}^2} = 0 \quad \text{ for each } j = 2, 3, \ldots, J, \text{ and} \\
\displaystyle\frac{\partial \mathscr{L}}{\partial \zeta} &= \frac{V_{\level_1}}{N_{\level_1}} +\sum_{j=2}^J \frac{V_{\level_j, \level_{j-1}}}{N_{\level_j}} -  \frac{\varepsilon^2}{2} = 0.
\end{aligned}
\end{cases}
\end{equation}
The solution of this system of equations is
\begin{align}
N_{\level_1} &= \frac{2}{{\varepsilon^2}} \sqrt{\frac{V_{\level_1}}{C_{\level_1}}} \left(\sqrt{V_{\level_1}C_{\level_1}} + \sum_{j=2}^J \sqrt{V_{\level_j, \level_{j-1}} C_{\level_j, \level_{j-1}}}\right), \text{ and}\label{eq:Nopt1}\\
N_{\level_j} &= \frac{2}{{\varepsilon^2}} \sqrt{\frac{V_{\level_j, \level_{j-1}}}{C_{\level_j, \level_{j-1}}}} \left(\sqrt{V_{\level_1}C_{\level_1}} + \sum_{j=2}^J \sqrt{V_{\level_j, \level_{j-1}} C_{\level_j, \level_{j-1}}}\right) \text{ for } j = 2, 3, \ldots, J.\label{eq:Noptj}
\end{align}
In an actual implementation of the ML-KDMC estimator, these values must be rounded up to the nearest integer to enforce an integer number of samples, and sample variances and cost estimates can replace the quantities $V_{\level_1}$, $V_{\level_j, \level_{j-1}}$, $C_{\level_1}$ and $C_{\level_j, \level_{j-1}}$, see, e.g.,~\cite{giles15}.

Substituting the optimal values for $N_{\level_1}$, $N_{\level_2}$, $\ldots$, $N_{\level_J}$ from~\eqref{eq:Noptj} into the total cost of the MLMC estimator, we find that
\begin{equation}\label{eq:total_cost}
\cost{\VAqoiest_{\calL}} = \frac{2}{{\varepsilon^2}} \left(\sqrt{V_{\level_1} C_{\level_1}} + \sum_{j=2}^J \sqrt{V_{\level_j, \level_{j-1}} C_{\level_j, \level_{j-1}}}\right)^2.
\end{equation}
We will use this expression for the total cost of the estimator in the next section, when computing the optimal selection of levels $\calL$ in~\secref{subsec:levelselect} below.

We remark that, for the KDMC scheme outlined in~\secref{sec:ml_sim}, the behanviour of the variances  $V_{\level_j, \level_{j-1}}$ and costs $C_{\level_j, \level_{j-1}}$ is highly nontrivial, and different from the usual monotone behanviour in the case of SDEs, as shown in~\cite{giles08}, or in the case of PDEs with random coefficients, as shown in, e.g.,~\cite{cliffe11}. See~\figref{fig:V-B1} and the discussion in~\secref{subsec:numlevelselect} below for details. This nontrivial behanviour poses an additional difficulty when selecting the level set $\calL$ below. Furthermore, the standard theoretical convergence results for MLMC, as presented in, e.g.,~\cite{giles08}, cannot be used in our case. However, our numerical results in~\secref{sec:numres} illustrate that our method achieves the usual cost complexity rate $\order{\varepsilon^{-2}}$, where $\varepsilon$ is the tolerance on the RMSE.

\subsection{Level selection}\label{subsec:levelselect}

In this section, we are looking for the subset of levels $\calL^\star \in S_\calL$ that yields the MLMC estimator with smallest possible cost, while ensuring that the bias constraint, i.e., equation~\eqref{eq:bias_constraint} is satisfied. The latter constraint can be satisfied by choosing a suitable most accurate level $L$. This value $L$ will be larger for smaller tolerances $\varepsilon^2$ imposed on the MSE, i.e., equation~\eqref{eq:mse}. Since, by construction, any feasible subset $\calL$ contains the most accurate level $L$, we choose $L$ such that constraint~\eqref{eq:bias_constraint} is satisfied. 

Using the expression for the total cost of the MLMC estimator from equation~\eqref{eq:total_cost}, the optimal subset of levels $\calL^\star$ is
\begin{equation}\label{eq:level_selection}
\calL^\star = \argmin_{\calL\in S_\calL} \left( \sqrt{V_{\level_1} C_{\level_1}} + \sum_{j=2}^J \sqrt{V_{\level_j, \level_{j-1}} C_{\level_j, \level_{j-1}}} \right).
\end{equation}
Thus, $\calL^\star$ is the solution of a combinatorial optimization problem. This problem could be solved using a brute-force approach, where we compute the value of the cost function for every feasible subset $\calL \in S_\calL$, for a given finest level $L$. However, this approach quickly becomes intractable, even for moderate values of $L$. Accordingly, let us introduce the dummy variables
\begin{align}
u_i &\coloneqq \begin{cases}
1 &\text{if } \level_1 = i\\
0 &\text{otherwise}
\end{cases}
&&\text{ for } i = 0, 1, \ldots, L, \text{ and} \\
w_{i,j} &\coloneqq \begin{cases}
1 &\text{if } \level_i \in \calL \wedge \level_j \in \calL \\
0 &\text{otherwise}
\end{cases}
&&\text{ for } i = 1, 2, \ldots, L \text{ and } j = 0, 1, \ldots, i - 1.
\end{align}
The unconstrained combinatorial optimization problem in equation~\eqref{eq:level_selection} can be written as a constrained integer linear programming problem
\begingroup
\allowdisplaybreaks
\begin{align}\label{eq:integerLP}
\min_{u_i, w_{i, j}} &\left( \sum_{i=0}^L \left(u_i \sqrt{V_{i} C_{i}}\right) + \sum_{i=1}^L \sum_{j=0}^{i-1} \left( w_{i, j} \sqrt{V_{i, j} C_{i, j}} \right) \right), \\
\text{subject to  } &u_i \in \{0, 1\}   &&i = 0, 1, \ldots, L,\\
&w_{i, j} \in \{0, 1\} &&i = 1, 2, \ldots, L, \; j = 0, 1, \ldots, i - 1\\
&\!\!\!\sum_{i=j+1}^L w_{i, j} \leq 1 &&j=0, 1, \ldots, L-1,\\
&\sum_{j=0}^{i-1} \;w_{i, j} \leq 1 &&i=1, 2, \ldots, L-1,\\
&\sum_{j=0}^{L-1} \;w_{L, j} = 1 &&\\
&\sum_{j=0}^{k-1} \;w_{k, j} - \sum_{i=k+1}^L w_{i, k} = 0 &&k=1, 2, \ldots, L-1,\\
&\sum_{i=0}^{L} \;u_i = 1 &&i=0, 1, \ldots, L.
\end{align}
\endgroup

The optimal set of levels $\calL^\star$ then simply consists of all levels $i$ where $u_i=1$ or $w_{i,j}=1$. Dedicated methods exist for solving~\eqref{eq:integerLP}, including \emph{branch-and-bound} methods, see, e.g.,~\cite{applegate06}, and \emph{cutting plane} methods, see, e.g.,~\cite{padberg91}.  However, this approach suffers from a major drawback: it assumes that the $L(L-1)/2$ values for the variances $V_{i, j},  i = 1, 2, \ldots, L$ and $j = 0, 1, \ldots, i - 1$, are available. We want to avoid estimating these quantities directly, and want to devise a strategy that uses approximate values for $V_{i, j}$ that can be obtained from a pilot run with a limited number of samples. By definition, the variance $V_{i, j}$ can be expressed as
\begin{equation}
V_{i, j} = V_{i} + V_{j} - 2 \rho_{i, j} \sqrt{V_{i}, V_{j}},
\end{equation}
where $\rho_{i, j}$ is the correlation between $\VAqoi_i$ and $\VAqoi_j$. Comparing $V_{L, i}$ and $V_{L, j}$ with $V_{i, j}$ and using the above expression, we find that
\begin{equation}
\rho_{L, i} \rho_{L, j} - \sqrt{\rho_{L, i}^2 \rho_{L, j}^2 + 1 - \rho_{L, i}^2 - \rho_{L, j}^2} \leq \rho_{i, j} \leq \rho_{L, i} \rho_{L, j} + \sqrt{\rho_{L, i}^2 \rho_{L, j}^2 + 1 - \rho_{L, i}^2 - \rho_{L, j}^2},
\end{equation}
where $\rho_{L, i}$ is the correlation coefficient between $\VAqoi_L$ and $\VAqoi_i$, and $\rho_{L, j}$ is the correlation coefficient between $\VAqoi_L$ and $\VAqoi_j$. Assuming that $\rho_{i, j}$ is the geometric mean of both extremes, see~\cite{giles08}, we find that
\begin{equation}
\rho_{i, j}^2 + 1 \approx \rho_{L, i}^2 + \rho_{L, j}^2,
\end{equation}
and hence
\begin{equation}\label{eq:V_approx}
V_{i, j} \approx V_{i} + V_{j} - 2\sqrt{(\rho_{L, i}^2 + \rho_{L, j}^2 - 1)V_{i}V_{j}}.
\end{equation}
The latter expression can be obtained from the $L+1$ values $V_{i}, \, i=0, 1, \ldots, L$ and the $L$ values $\rho_{L, i},\, i = 0, 1, \ldots, L-1$. In our numerical experiments presented below, we use sample variances and sample correlation coefficients to approximate $V_{i}$ and $\rho_{L, i}$,  that are extracted from a pilot run with a limited number of samples on each level. From this set of samples, we can also extract actual run times that can replace the cost estimates $C_{i}$ and $C_{i, j}$ in equation~\eqref{eq:integerLP}. In~\secref{subsec:numlevelselect}, we will use the strategy outlined in this section to derive some general guidelines for level selection, that work for a large range of plasma background configurations.

\section{Numerical results}\label{sec:numres}

\begin{figure}
\centering
	\tikzexternalenable
	\tikzsetnextfilename{backgrounds}
\begin{tikzpicture}
	\begin{groupplot}[%
		group style={
			group size=2 by 1,
			 vertical sep=1.5em,
			 horizontal sep=1.5em
		},
		default axis
	]
	\nextgroupplot[%
		xmin=-4, xmax=4, ymin=98, ymax=104, ytick={98, 100, 102, 104}, xtick={-4,-2,...,4}, xticklabels={$-1$, $0$, $1$, $2$, $3$}, yticklabels={, $b$, $b(a+1)$, $b(2a+1)$}
	]
	\def\a{1}
	\def\b{100}
	\def\s{1}
	\addplot[domain=-4:0, thick, line cap=round, color={red}] {-\a*(x-\s) + \b - \a*\s};
	\addplot[domain=0:4, thick, line cap=round, color={red}] {\a*(x-\s) + \b + \a*\s};
	\coordinate (t1) at (axis description cs:.5,1.2);
	\pgfplotsset{
		after end axis/.append code={%
			\node [] at (t1) {collision rate $R_1(x)$};
		}
	}
	\nextgroupplot[%
		xmin=-4, xmax=4, ymin=98, ymax=104, ytick={98, 100, 102, 104}, xtick={-4,-2,...,4}, xticklabels={$-1$, $0$, $1$, $2$, $3$}, yticklabels=\empty
	]
	\def\a{1}
	\def\b{100}
	\def\s{1}
	\addplot[domain=-4:0, thick, line cap=round, color={red}] {\b};
	\addplot[domain=0:4, thick, line cap=round, color={red}] {\a*(x-\s) + \b + \a*\s};
	\coordinate (t2) at (axis description cs:.5,1.2);
	\pgfplotsset{
		after end axis/.append code={%
			\node [] at (t2) {collision rate $R_2(x)$};
		}
	}
\end{groupplot}
\end{tikzpicture}
	\tikzexternaldisable

\caption{Collision rate $R_1(x)$ for the symmetric background $\B{1}$ (\emph{left}) and collision rate $R_2(x)$ for the asymmetric background $\B{2}$ (\emph{right}).}
\label{fig:backgrounds}
\end{figure}
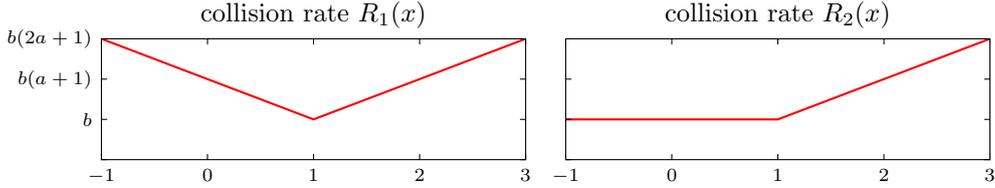

In this section, we apply the ML-KDMC method to compute the expected value of the particle position for a one-dimensional test-case inspired by~\cite{dimarco18}. We set up a simulation for $t \in [0, \VAendtime]$ with end time $\VAendtime=1$. Particles are released from their initial position $x(0)=1$, and collide with two different backgrounds $\B{1}\coloneqq\{R_1(x), \VAbackgroundspeedmean, \VAbackgroundspeedstdv\}$ and $\B{2}\coloneqq\{R_2(x), \VAbackgroundspeedmean, \VAbackgroundspeedstdv\}$ with $\VAbackgroundspeedmean=0$, $\VAbackgroundspeedstdv=1$ and collision rates
\begin{equation}
R_1(x) = \begin{cases}
-b (a(x - 1) - 1) & x \le 1 \\
\hphantom{-}b(a(x - 1) + 1) & x > 1
\end{cases} \quad \text{and} \quad
R_2(x) = \begin{cases}
\hphantom{-}b & x \le 1 \\
\hphantom{-}b(a(x - 1) + 1) & x > 1
\end{cases}.
\end{equation}
These collision rates are shown in~\figref{fig:backgrounds}. For obvious reasons, we call $\B{1}$ a \emph{symmetric} background, and $\B{2}$ an \emph{asymmetric} background. We refer to the background with $a=0$ as the \emph{homogeneous} case, since, in that case, the collision rate is constant, i.e., we have $R_1(x)=R_2(x)=b$. A background with $a\gg0$ is referred to as the \emph{heterogeneous} case. We choose the background parameters as  $a \in \{0, 0.1, 0.2, 0.5, 1, 2, 5, 10, 20, 50, 100\}$ and $b \in\{1, 10, 100, 1\,000, 10\,000, 100\,000\}$.

\setlength{\figurewidth}{.45\textwidth}
\setlength{\figureheight}{.2\textwidth}

Our main numerical results are divided into three subsections. First, in~\secref{subsec:consistency}, we will check the consistency of the random number mapping and aggregation $\varphi_{\level_j}^{\level_{j-1}} : \VArandparam_{\level_j}^{(\VAparticleno)} \mapsto \tilde\VArandparam_{\level_{j-1}}^{(\VAparticleno)}$ from~\secref{sec:ml_cor} numerically. That is, we will verify if the aggregated random numbers $\tilde\VArandparam_{\level_{j-1}}^{(\VAparticleno)}$ for the coarse particle path satisfy the required distributions. Next, in~\secref{subsec:numlevelselect}, we study the level selection strategy from~\secref{subsec:levelselect}, and devise a heuristic algorithm for level selection in the ML-KDMC method. Finally, in~\secref{subsec:performance}, we compare the efficiency of the new ML-KDMC scheme with the standard, single-level KDMC scheme.

In all our numerical experiments below, we used our implementation of the ML-KDMC scheme, available online at~\url{https://numa.cs.kuleuven.be/research_private/pieterjan.robbe/apmlmc/ML-APMC.jl.tar.gz}, and the MLMC code \texttt{Mul\-ti\-le\-vel\-Es\-ti\-ma\-tors}, available online at~\url{https://numa.cs.kuleuven.be/research_private/pieterjan.robbe/apmlmc/MultilevelEstimators.jl.tar.gz}.

\subsection{Coarse particle path consistency}\label{subsec:consistency}

In this section, we numerically verify the consistency of the random numbers used for the coarse correlated particle path. That is, we will assert whether the random numbers $\tilde\VArandparam_{\level_{j-1}}^{(\VAparticleno)}$, defined by the mapping $\varphi_{\level_j}^{\level_{j-1}}$ from~\secref{sec:ml_cor}, follow the required distributions. We will use the Anderson--Darling (AD) hypothesis test to test whether a given set of samples is drawn from a certain probability distribution. Let $\Phi(x)$ denote the cumulative distribution function (CDF) of the proposed distribution, and let $\Phi_N(x)$ denote the empirical CDF obtained from $N$ samples $\{\xi^{(\VAparticleno)}\}_{n=1}^N$ of the random variable $\xi$. Recall that the latter is computed as
\begin{equation*}
\Phi_N(x) \coloneqq \frac{1}{N} \sum_{n=1}^N \mathbbm{1}_{\xi^{(\VAparticleno)} \le x},
\end{equation*}
where $\mathbbm{1}_{\xi^{(\VAparticleno)} \le x}$ is the indicator function for event $\xi^{(\VAparticleno)} \le x$. The AD test computes the metric
\begin{equation}\label{eq:anderson}
A^2 = N \int_{-\infty}^{\infty} w(x) (\Phi_N(x) - \Phi(x))^2 \mathrm{d} \Phi(x) \quad \text{where} \quad w(x) = \frac{1}{\Phi(x)(1 - \Phi(x))}.
\end{equation}
A larger value of the distance $A^2$ means that it is less likely that the samples are coming from the proposed distribution with CDF $\Phi$. The inference problem can be solved using the null hypothesis
\begin{equation}
H_0 : \xi \text{ follows a distribution with CDF } \Phi(x).
\end{equation}
To assert the validity of the null hypothesis, we compute the so-called $p$-value, i.e., the probability that, under the proposed distribution in the null hypothesis, the value of $A^2$ is at least as large as the value of $A^2$ that was computed from the available samples. The null hypothesis $H_0$ is rejected only when this $p$-value is below a certain threshold, say 1\%.

Let us repeat these steps for the three random variables $\tilde\VArnvel$, $\tilde\VArnexp$ and $\tilde\VArndif$, where we dropped the subscript $(\level-1, \VAeventno)$ for convenience. The corresponding null hypotheses are
\begin{equation}\label{eq:null_hypothesis}\begin{aligned}
H_0^{\tilde\VArnvel} &: \tilde\VArnvel \text{ follows a normal distribution with CDF } \Phi^\calN(x)\text{, i.e., } \tilde\VArnvel \sim \calN(0, 1), \\
H_0^{\tilde\VArnexp} &: \tilde\VArnexp \text{ follows an exponential distribution with CDF } \Phi^\calE(x)\text{, i.e., } \tilde\VArnexp \sim \calE(1), \text{ and} \\
H_0^{\tilde\VArndif} &: \tilde\VArndif \text{ follows a normal distribution with CDF } \Phi^\calN(x)\text{, i.e., } \tilde\VArndif \sim \calN(0, 1).
\end{aligned}
\end{equation}
\tabref{tab:hypothesis} shows the AD distance $A^2$ and the $p$-value for background $\B{1}$ with $a=10$ and $b=100$ for various (fine) levels $\level$, based on $N=1\,000$ samples. The time step is given by $\VAtimestep_\level\coloneqq\VAendtime/2^\level$ for the fine particle, and by $\VAtimestep_{\level-1}\coloneqq\VAfactor \VAtimestep_\level$ for the coarse particle. The entry $\level=5$ and $\VAfactor=2$, for example, means that we look for the random numbers $\tilde\VArandparam_4$ obtained from a coarsening of the random numbers $\VArandparam_5$ on level $\level=5$, with coarsening factor $M=2$. All samples are  obtained by running repeated particle simulations and recording all values for the (coarse particle) random numbers $\tilde\VArnexp_{\level-1, \VAeventno}$, $\tilde\VArnvel_{\level-1, \VAeventno}$ and $\tilde\VArndif_{\level-1, \VAeventno}$ for each collision $\VAeventno=0,\ldots,\tilde\VAncollisions_{\level_{j-1}}^{(\VAparticleno)}$, until $N=1\,000$ realizations are available. We observe that in all cases, we fail to reject the null hypothesis, i.e., the $p$-value is above 0.01 (1\%). Hence, we accept the hypotheses in equation~\eqref{eq:null_hypothesis}, and find that the coarse aggregated random numbers satisfy the required distributions. 

\begin{figure}[b]
\centering
	\tikzexternalenable
	\tikzsetnextfilename{plot_cdfs}
\pgfmathtruncatemacro{\a}{10}
\pgfmathtruncatemacro{\b}{1000}
\pgfmathtruncatemacro{\id}{1}
\pgfmathtruncatemacro\level{5}
\pgfmathtruncatemacro\k{1}
\pgfmathtruncatemacro{\mykval}{2^\k}
\setlength{\figurewidth}{.33\textwidth}
\setlength{\figureheight}{.33\textwidth}

\begin{tikzpicture}
	\begin{groupplot}[%
		group style={
			group size=3 by 1,
			 horizontal sep=3em,
		},
		default axis
	]
	\nextgroupplot[%
		xmin=-3, xmax=3, ymin=0, ymax=1, ylabel={$\Phi(x)$}, xlabel={$x$}, xtick={-3, -2, ..., 3}, title={$\tilde\nu$}
	]
	\addplot [line cap=round, thick, red] table[header=false, x index=0, y index=1] {../JuliaCode/data/normal_cdf}; \addlegendentry{exact}
	\addplot [line cap=round] table[header=false, x index=0, y index=2] {data/emperical_cdf_B1_a_\a_b_\b_ell_\level_k_\k}; \addlegendentry{empirical}
	\nextgroupplot[%
		xmin=0, xmax=6, ymin=0, ymax=1, xlabel={$x$}, xtick={0, 1, ..., 6}, yticklabels=\empty, title={$\tilde\epsilon$}
	]
	\addplot [line cap=round, thick, green!70!black] table[header=false, x index=0, y index=1] {../JuliaCode/data/exponential_cdf}; \addlegendentry{exact}
	\addplot [line cap=round] table[header=false, x index=1, y index=3] {data/emperical_cdf_B1_a_\a_b_\b_ell_\level_k_\k}; \addlegendentry{empirical}
	\nextgroupplot[%
		xmin=-3, xmax=3, ymin=0, ymax=1, xlabel={$x$}, xtick={-3, -2, ..., 3}, yticklabels=\empty, title={$\tilde\chi$}
	]
	\addplot [line cap=round, thick, azure] table[header=false, x index=0, y index=1] {../JuliaCode/data/normal_cdf}; \addlegendentry{exact}
	\addplot [line cap=round] table[header=false, x index=0, y index=4] {data/emperical_cdf_B1_a_\a_b_\b_ell_\level_k_\k}; \addlegendentry{empirical}
\end{groupplot}
\end{tikzpicture}
	\tikzexternaldisable

\caption{A visual comparison of the exact CDF and its empirical counterpart for the random variables $\tilde\VArnvel$, $\tilde\VArnexp$ and $\tilde\VArndif$ for background $\B{1}$ with $a=10$, $b=1\,000$, $\ell=5$ and $K=2$.}
\label{fig:cdfs}
\end{figure}
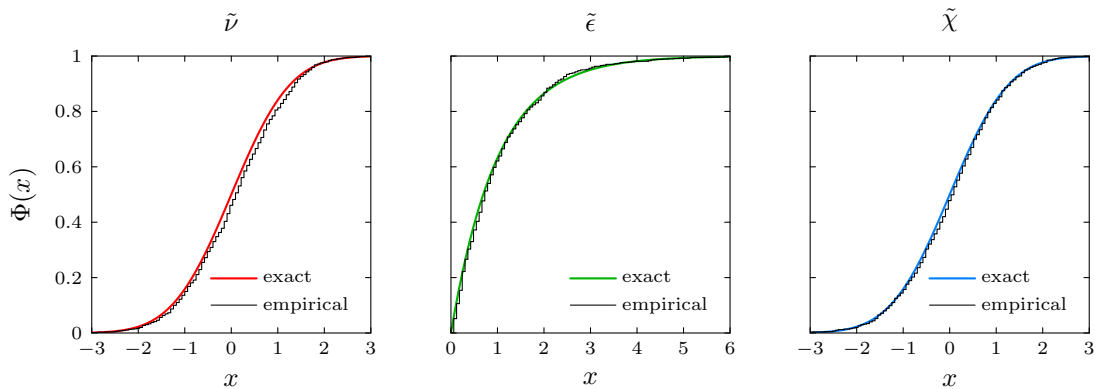

\begin{table}\centering\small
\begin{tabular}{llcccccccccccc}\toprule
&& \multicolumn{2}{c}{$\VAfactor = 2$} & \multicolumn{2}{c}{$\VAfactor = 4$} & \multicolumn{2}{c}{$\VAfactor = 8$} & \multicolumn{2}{c}{$\VAfactor = 16$} & \multicolumn{2}{c}{$\VAfactor = 32$} & \multicolumn{2}{c}{$\VAfactor = 64$} \\\midrule
&& $A^2$ & $p$ & $A^2$ & $p$ & $A^2$ & $p$ & $A^2$ & $p$ & $A^2$ & $p$ & $A^2$ & $p$ \\ \midrule
\multirow{3}{*}{\rotatebox[origin=c]{90}{$\ell=1$}}& $\tilde\nu$ & 0.344 & \cellcolor{yescolor}0.902& & & & & & & & & & \\
& $\tilde\epsilon$ & 0.622 & \cellcolor{yescolor}0.627& & & & & & & & & & \\
& $\tilde\chi$ & 0.207 & \cellcolor{yescolor}0.988& & & & & & & & & & \\
\midrule\multirow{3}{*}{\rotatebox[origin=c]{90}{$\ell=2$}}& $\tilde\nu$ & 0.171 & \cellcolor{yescolor}0.996& 0.656 & \cellcolor{yescolor}0.597& & & & & & & & \\
& $\tilde\epsilon$ & 0.352 & \cellcolor{yescolor}0.895& 0.397 & \cellcolor{yescolor}0.852& & & & & & & & \\
& $\tilde\chi$ & 0.385 & \cellcolor{yescolor}0.863& 0.691 & \cellcolor{yescolor}0.566& & & & & & & & \\
\midrule\multirow{3}{*}{\rotatebox[origin=c]{90}{$\ell=3$}}& $\tilde\nu$ & 0.286 & \cellcolor{yescolor}0.948& 0.422 & \cellcolor{yescolor}0.827& 0.265 & \cellcolor{yescolor}0.962& & & & & & \\
& $\tilde\epsilon$ & 0.553 & \cellcolor{yescolor}0.693& 0.574 & \cellcolor{yescolor}0.673& 0.427 & \cellcolor{yescolor}0.822& & & & & & \\
& $\tilde\chi$ & 0.165 & \cellcolor{yescolor}0.997& 0.364 & \cellcolor{yescolor}0.884& 0.554 & \cellcolor{yescolor}0.692& & & & & & \\
\midrule\multirow{3}{*}{\rotatebox[origin=c]{90}{$\ell=4$}}& $\tilde\nu$ & 0.329 & \cellcolor{yescolor}0.915& 0.626 & \cellcolor{yescolor}0.624& 0.440 & \cellcolor{yescolor}0.809& 0.602 & \cellcolor{yescolor}0.646& & & & \\
& $\tilde\epsilon$ & 0.415 & \cellcolor{yescolor}0.834& 0.296 & \cellcolor{yescolor}0.941& 1.088 & \cellcolor{yescolor}0.314& 0.483 & \cellcolor{yescolor}0.764& & & & \\
& $\tilde\chi$ & 0.362 & \cellcolor{yescolor}0.885& 0.242 & \cellcolor{yescolor}0.974& 0.440 & \cellcolor{yescolor}0.808& 0.310 & \cellcolor{yescolor}0.930& & & & \\
\midrule\multirow{3}{*}{\rotatebox[origin=c]{90}{$\ell=5$}}& $\tilde\nu$ & 0.261 & \cellcolor{yescolor}0.964& 0.258 & \cellcolor{yescolor}0.966& 0.693 & \cellcolor{yescolor}0.565& 0.629 & \cellcolor{yescolor}0.621& 0.251 & \cellcolor{yescolor}0.970& & \\
& $\tilde\epsilon$ & 0.131 & \cellcolor{yescolor}1.000& 0.826 & \cellcolor{yescolor}0.463& 0.196 & \cellcolor{yescolor}0.991& 0.574 & \cellcolor{yescolor}0.673& 0.512 & \cellcolor{yescolor}0.735& & \\
& $\tilde\chi$ & 0.277 & \cellcolor{yescolor}0.954& 0.402 & \cellcolor{yescolor}0.847& 0.166 & \cellcolor{yescolor}0.997& 0.425 & \cellcolor{yescolor}0.824& 0.611 & \cellcolor{yescolor}0.637& & \\
\midrule\multirow{3}{*}{\rotatebox[origin=c]{90}{$\ell=6$}}& $\tilde\nu$ & 0.550 & \cellcolor{yescolor}0.696& 0.286 & \cellcolor{yescolor}0.948& 0.473 & \cellcolor{yescolor}0.774& 0.433 & \cellcolor{yescolor}0.816& 0.414 & \cellcolor{yescolor}0.835& 0.416 & \cellcolor{yescolor}0.833\\
& $\tilde\epsilon$ & 1.255 & \cellcolor{yescolor}0.248& 0.577 & \cellcolor{yescolor}0.670& 0.344 & \cellcolor{yescolor}0.902& 0.478 & \cellcolor{yescolor}0.769& 0.610 & \cellcolor{yescolor}0.639& 0.243 & \cellcolor{yescolor}0.974\\
& $\tilde\chi$ & 0.290 & \cellcolor{yescolor}0.945& 0.813 & \cellcolor{yescolor}0.472& 0.301 & \cellcolor{yescolor}0.938& 0.260 & \cellcolor{yescolor}0.965& 0.671 & \cellcolor{yescolor}0.583& 0.271 & \cellcolor{yescolor}0.958\\
\midrule\multirow{3}{*}{\rotatebox[origin=c]{90}{$\ell=7$}}& $\tilde\nu$ & 0.365 & \cellcolor{yescolor}0.882& 0.348 & \cellcolor{yescolor}0.898& 0.382 & \cellcolor{yescolor}0.866& 0.340 & \cellcolor{yescolor}0.905& 0.135 & \cellcolor{yescolor}0.999& 0.316 & \cellcolor{yescolor}0.926\\
& $\tilde\epsilon$ & 0.648 & \cellcolor{yescolor}0.604& 0.254 & \cellcolor{yescolor}0.968& 1.299 & \cellcolor{yescolor}0.233& 0.426 & \cellcolor{yescolor}0.822& 0.282 & \cellcolor{yescolor}0.951& 0.413 & \cellcolor{yescolor}0.836\\
& $\tilde\chi$ & 0.463 & \cellcolor{yescolor}0.785& 0.184 & \cellcolor{yescolor}0.994& 0.312 & \cellcolor{yescolor}0.929& 0.174 & \cellcolor{yescolor}0.996& 0.233 & \cellcolor{yescolor}0.979& 0.254 & \cellcolor{yescolor}0.968\\
\midrule\multirow{3}{*}{\rotatebox[origin=c]{90}{$\ell=8$}}& $\tilde\nu$ & 0.682 & \cellcolor{yescolor}0.574& 0.252 & \cellcolor{yescolor}0.969& 0.416 & \cellcolor{yescolor}0.833& 0.250 & \cellcolor{yescolor}0.970& 0.641 & \cellcolor{yescolor}0.610& 0.103 & \cellcolor{yescolor}1.000\\
& $\tilde\epsilon$ & 0.561 & \cellcolor{yescolor}0.685& 1.043 & \cellcolor{yescolor}0.335& 0.533 & \cellcolor{yescolor}0.713& 0.336 & \cellcolor{yescolor}0.909& 1.183 & \cellcolor{yescolor}0.274& 0.566 & \cellcolor{yescolor}0.680\\
& $\tilde\chi$ & 0.179 & \cellcolor{yescolor}0.995& 0.707 & \cellcolor{yescolor}0.553& 0.325 & \cellcolor{yescolor}0.918& 0.511 & \cellcolor{yescolor}0.736& 0.253 & \cellcolor{yescolor}0.969& 0.295 & \cellcolor{yescolor}0.942\\
\midrule\multirow{3}{*}{\rotatebox[origin=c]{90}{$\ell=9$}}& $\tilde\nu$ & 0.434 & \cellcolor{yescolor}0.814& 0.234 & \cellcolor{yescolor}0.978& 0.421 & \cellcolor{yescolor}0.828& 0.293 & \cellcolor{yescolor}0.943& 0.554 & \cellcolor{yescolor}0.692& 0.288 & \cellcolor{yescolor}0.947\\
& $\tilde\epsilon$ & 0.748 & \cellcolor{yescolor}0.520& 0.560 & \cellcolor{yescolor}0.687& 0.305 & \cellcolor{yescolor}0.935& 0.200 & \cellcolor{yescolor}0.991& 0.545 & \cellcolor{yescolor}0.702& 0.408 & \cellcolor{yescolor}0.841\\
& $\tilde\chi$ & 0.702 & \cellcolor{yescolor}0.557& 0.530 & \cellcolor{yescolor}0.717& 0.478 & \cellcolor{yescolor}0.770& 0.736 & \cellcolor{yescolor}0.529& 0.351 & \cellcolor{yescolor}0.895& 0.458 & \cellcolor{yescolor}0.790\\
\midrule\multirow{3}{*}{\rotatebox[origin=c]{90}{$\ell=10$}}& $\tilde\nu$ & 0.312 & \cellcolor{yescolor}0.929& 0.563 & \cellcolor{yescolor}0.684& 0.215 & \cellcolor{yescolor}0.986& 0.628 & \cellcolor{yescolor}0.622& 0.644 & \cellcolor{yescolor}0.607& 0.633 & \cellcolor{yescolor}0.618\\
& $\tilde\epsilon$ & 0.359 & \cellcolor{yescolor}0.888& 0.594 & \cellcolor{yescolor}0.653& 1.046 & \cellcolor{yescolor}0.334& 0.585 & \cellcolor{yescolor}0.662& 0.457 & \cellcolor{yescolor}0.791& 0.187 & \cellcolor{yescolor}0.994\\
& $\tilde\chi$ & 0.234 & \cellcolor{yescolor}0.978& 0.376 & \cellcolor{yescolor}0.873& 0.194 & \cellcolor{yescolor}0.992& 0.301 & \cellcolor{yescolor}0.938& 0.294 & \cellcolor{yescolor}0.943& 0.343 & \cellcolor{yescolor}0.903\\
\midrule\multirow{3}{*}{\rotatebox[origin=c]{90}{$\ell=11$}}& $\tilde\nu$ & 0.824 & \cellcolor{yescolor}0.464& 0.585 & \cellcolor{yescolor}0.662& 0.382 & \cellcolor{yescolor}0.866& 0.623 & \cellcolor{yescolor}0.626& 0.252 & \cellcolor{yescolor}0.969& 0.234 & \cellcolor{yescolor}0.978\\
& $\tilde\epsilon$ & 0.637 & \cellcolor{yescolor}0.614& 0.438 & \cellcolor{yescolor}0.810& 0.332 & \cellcolor{yescolor}0.912& 0.567 & \cellcolor{yescolor}0.680& 1.101 & \cellcolor{yescolor}0.309& 0.307 & \cellcolor{yescolor}0.932\\
& $\tilde\chi$ & 0.350 & \cellcolor{yescolor}0.896& 0.291 & \cellcolor{yescolor}0.945& 0.331 & \cellcolor{yescolor}0.914& 0.575 & \cellcolor{yescolor}0.672& 0.237 & \cellcolor{yescolor}0.977& 0.575 & \cellcolor{yescolor}0.672\\
\midrule\multirow{3}{*}{\rotatebox[origin=c]{90}{$\ell=12$}}& $\tilde\nu$ & 0.258 & \cellcolor{yescolor}0.966& 0.717 & \cellcolor{yescolor}0.545& 0.206 & \cellcolor{yescolor}0.989& 0.483 & \cellcolor{yescolor}0.764& 0.412 & \cellcolor{yescolor}0.836& 0.250 & \cellcolor{yescolor}0.970\\
& $\tilde\epsilon$ & 0.595 & \cellcolor{yescolor}0.653& 0.293 & \cellcolor{yescolor}0.943& 0.560 & \cellcolor{yescolor}0.686& 0.357 & \cellcolor{yescolor}0.890& 0.756 & \cellcolor{yescolor}0.514& 0.273 & \cellcolor{yescolor}0.957\\
& $\tilde\chi$ & 0.265 & \cellcolor{yescolor}0.962& 0.481 & \cellcolor{yescolor}0.767& 0.252 & \cellcolor{yescolor}0.969& 0.243 & \cellcolor{yescolor}0.974& 0.300 & \cellcolor{yescolor}0.938& 0.748 & \cellcolor{yescolor}0.520\\
\midrule\multirow{3}{*}{\rotatebox[origin=c]{90}{$\ell=13$}}& $\tilde\nu$ & 1.012 & \cellcolor{yescolor}0.351& 0.257 & \cellcolor{yescolor}0.966& 0.320 & \cellcolor{yescolor}0.922& 0.404 & \cellcolor{yescolor}0.844& 0.301 & \cellcolor{yescolor}0.937& 0.374 & \cellcolor{yescolor}0.874\\
& $\tilde\epsilon$ & 0.790 & \cellcolor{yescolor}0.488& 0.709 & \cellcolor{yescolor}0.552& 1.787 & \cellcolor{yescolor}0.121& 0.418 & \cellcolor{yescolor}0.831& 0.477 & \cellcolor{yescolor}0.770& 0.458 & \cellcolor{yescolor}0.789\\
& $\tilde\chi$ & 0.185 & \cellcolor{yescolor}0.994& 0.697 & \cellcolor{yescolor}0.561& 0.110 & \cellcolor{yescolor}1.000& 0.142 & \cellcolor{yescolor}0.999& 0.281 & \cellcolor{yescolor}0.951& 0.412 & \cellcolor{yescolor}0.836\\
\midrule\multirow{3}{*}{\rotatebox[origin=c]{90}{$\ell=14$}}& $\tilde\nu$ & 0.632 & \cellcolor{yescolor}0.618& 0.389 & \cellcolor{yescolor}0.860& 0.355 & \cellcolor{yescolor}0.892& 0.584 & \cellcolor{yescolor}0.663& 0.265 & \cellcolor{yescolor}0.962& 0.584 & \cellcolor{yescolor}0.664\\
& $\tilde\epsilon$ & 0.253 & \cellcolor{yescolor}0.969& 0.806 & \cellcolor{yescolor}0.477& 0.266 & \cellcolor{yescolor}0.961& 0.229 & \cellcolor{yescolor}0.980& 0.662 & \cellcolor{yescolor}0.591& 0.379 & \cellcolor{yescolor}0.869\\
& $\tilde\chi$ & 0.557 & \cellcolor{yescolor}0.690& 0.126 & \cellcolor{yescolor}1.000& 0.296 & \cellcolor{yescolor}0.941& 0.360 & \cellcolor{yescolor}0.887& 0.235 & \cellcolor{yescolor}0.978& 0.150 & \cellcolor{yescolor}0.999\\
\midrule\multirow{3}{*}{\rotatebox[origin=c]{90}{$\ell=15$}}& $\tilde\nu$ & 0.962 & \cellcolor{yescolor}0.378& 0.499 & \cellcolor{yescolor}0.747& 0.375 & \cellcolor{yescolor}0.873& 0.232 & \cellcolor{yescolor}0.979& 0.269 & \cellcolor{yescolor}0.960& 0.255 & \cellcolor{yescolor}0.968\\
& $\tilde\epsilon$ & 0.333 & \cellcolor{yescolor}0.912& 0.561 & \cellcolor{yescolor}0.685& 0.421 & \cellcolor{yescolor}0.827& 0.356 & \cellcolor{yescolor}0.891& 0.608 & \cellcolor{yescolor}0.641& 0.258 & \cellcolor{yescolor}0.966\\
& $\tilde\chi$ & 0.390 & \cellcolor{yescolor}0.859& 0.415 & \cellcolor{yescolor}0.834& 0.192 & \cellcolor{yescolor}0.993& 0.370 & \cellcolor{yescolor}0.878& 0.340 & \cellcolor{yescolor}0.905& 0.507 & \cellcolor{yescolor}0.739\\
\bottomrule
\end{tabular}
\caption{Anderson--Darling distance $A^2$ and corresponding $p$-value for the aggregated coarse particle random numbers $\tilde\VArnexp$, $\tilde\VArnvel$ and $\tilde\VArndif$ for background $\B{1}$ with $a=10$ and $b=1\,000$ and various levels $\ell=1, 2, \ldots, 15$ and level multiplication factors $\VAfactor\in\{2, 4, 8, 16, 32, 64\}$ based on $N=1\,000$ samples. The random numbers do not follow the required distributions (i.e., the null hypothesis is rejected) when the $p$-value is below 0.05.}
\label{tab:hypothesis}
\end{table}

A visual comparison of the exact CDF and its empirical counterpart in case $\level=5$ and $\VAfactor=2$ is shown in~\figref{fig:cdfs} for all three random variables $\tilde\VArnvel$, $\tilde\VArnexp$ and $\tilde\VArndif$. Note that these results validate the telescopic sum in the multilevel estimator, i.e., equation~\eqref{eq:telsum}: since $\VArandparam_{\level_{j}}$ and $\tilde\VArandparam_{\level_{j-1}}$ follow the same distribution, we have that 
\begin{equation}
\E{\VAqoi_{\level_{j-1}}(\VArandparam_{\level_{j-1}})} = \E{\VAqoi_{\level_{j-1}}(\tilde\VArandparam_{\level_{j-1}})}.
\end{equation}
More results for background $\B{2}$ and for all other parameter combinations can be found online at \url{https://numa.cs.kuleuven.be/research_private/pieterjan.robbe/apmlmc/}.

\subsection{A heuristic level selection method}\label{subsec:numlevelselect}

Before we are able to use our ML-KDMC method, we must specify which set of levels $\calL$ should be used. This is an important decision, since the choice of the level set $\calL$ determines the efficiency of the multilevel estimator (it appears directly in the expression for the cost of the multilevel estimator, equation~\eqref{eq:total_cost}). Our strategy for selecting the optimal set of levels was outlined in~\secref{subsec:levelselect}. However, solving problem~\eqref{eq:integerLP} for all combinations of collision rate parameters $a$ and $b$ would be computationally infeasible. Instead, we compute the optimal set of levels $\calL^\star$ for a judiciously chosen set of parameter combinations with maximum level $L=22$, and hope to devise some general guidelines for level selection for all other parameter combinations. 

The effect of the level selection strategy can be visualized by inspecting the variances $V_\level$ and $V_{\level, \level-1}$, $\ell=0,1,\ldots,L$ and costs $C_\level$ and $C_{\level, \level-1}$, $\ell=0,1,\ldots,L$, see~Figures~\ref{fig:V-B1}--\ref{fig:C-B1}. These figures indicate a similar behanviour of the variance and cost of the multilevel difference for all parameter combinations, see the sketch in~\figref{fig:sketch}.

In the homogeneous case ($a=0$, \raisebox{.15\baselineskip}{\tikz{\draw[very thick, line cap=round, dashed] (0,0) -- (.5,0);}}), the variance of the multilevel difference $V_{\level, \level-1}$ increases with decreasing time step $\VAtimestep_\level$, until it reaches a maximum around $\VAtimestep_\level=1/R(x)$. For even smaller values of the time step $\VAtimestep_\level$, the variance decreases at a rate proportional to $\order{\VAtimestep_\level^3}$. As the heterogeneity increases ($a\uparrow$, \raisebox{.15\baselineskip}{\tikz{\draw[very thick, line cap=round] (0,0) -- (.5,0);}}), there is an additional error that decays as $\order{\VAtimestep_\level}$, which makes the variance of the multilevel difference decay to a local minimum before it increases again. We remark that the behanviour for $b=1$ and $b=10$ in Figure~\ref{fig:V-B1} is a degenerate case of the more general behanviour sketched above, where the local maximum is shifted to the right, and infeasible values of the time step size larger than the simulated time are required. The computational cost $C_{\level, \level-1}$ in~\figref{fig:C-B1} increases linearly with decreasing $\VAtimestep_\level$ until the local maximum in the variance decay curve is reached, after which the cost scales independent of the time step. This behanviour is consistent with the error analysis of the KDMC scheme in~\cite{mortier20}: for large values of the time step size $\VAtimestep_\level$, the scheme converges to the diffusive approximation of the Boltzmann-BGK equation, while for small values of $\VAtimestep_\level$, the scheme converges to the kinetic approximation (i.e., the direct simulation of each collision of the particle).

Based on these observations, we suggest the level selection strategy summarized in~\figref{fig:sketch-with-levels}. This level selection strategy is the result of applying the level selection method from~\secref{subsec:levelselect} to the simulation results shown in~Figures~\ref{fig:V-B1}--\ref{fig:C-B1}. Given the variance $V_\level$ for each $\level=0,1,\ldots,L$, and the variance of the multilevel difference $V_{\level, \level-1}$ for each $\level=1, 2, \ldots, L$, the set of levels $\calL$ can be selected heuristically using~\algref{alg:level_select}. This approach has several advantages compared to the level selection method from~\secref{subsec:levelselect}. First of all, there is no need to compute the solution of the integer linear programming problem from equation~\eqref{eq:integerLP} for every new background collision rate parameter combination. Furthermore, the required quantities $V_\level$ and $V_{\level, \level-1}$ can be approximated by the sample variance of a set of (cheap) warm-up samples for levels $\level$ where $\VAtimestep_\level \ge 1/R(x)$, and we can use the asymptotic complexity rate of $V_{\level, \level-1}$ for all other levels  $\level$ where $\VAtimestep_\level < 1/R(x)$. The warm-up samples for every level $\level\in\calL$ can be reused to hot start the multilevel estimator. Also, our approach is level-adaptive: increasing the maximum level parameter $L$ does not force us to recompute the level set $\calL$. Finally, the subsequent level sets $\calL$ are nested for different $L$, which means that no computational effort is wasted when $L$ is increased. This would be the case if a particular level $\level$ is part of the set $\calL$ for the maximum level parameter $L$, but is not a part of the set $\calL$ for the maximum level parameter $L+1$.

\setlength{\figurewidth}{.45\textwidth}
\setlength{\figureheight}{.45\textwidth}
\begin{figure}
\centering
	\tikzexternalenable
	\tikzsetnextfilename{V-B1}
%
%
\begin{tikzpicture}%
	\begin{groupplot}[%
		group style={
			group size=2 by 3,
			 vertical sep=2em,
			 horizontal sep=2em,
			 group name=my plots
		},
		loglog axis
	]
	\nextgroupplot[%
		xmin=1e-7, xmax=1e0, ymin=1e-18, ymax=1e2, ylabel={$\V{\;\cdot\;}$}, xticklabels=\empty
	]
	\foreach \a in {0, 0_1, 0_2, 0_5, 1, 2, 5, 10, 20, 50, 100}{
		\addplot+ [default line] table[header=false, x expr=1/2^\coordindex, y index=3] {data/analysis_B1_a_\a_b_1};
		\label{V-1-a\a}
	}
	\foreach \a in {0, 0_1, 0_2, 0_5, 1, 2, 5, 10, 20, 50, 100}{
		\addplot+ [hollow line] table[header=false, x expr=1/2^\coordindex, y index=2] {data/analysis_B1_a_\a_b_1};
	}
	\slopetriangle{.425}{.1}{.075}{3}
	\coordinate (t) at (axis description cs:.985,.03);
	\pgfplotsset{
		after end axis/.append code={%
			\node [anchor=south east] at (t) {$b = 1$};
		}
	}
	\nextgroupplot[%
		xmin=1e-7, xmax=1e0, ymin=1e-18, ymax=1e2, xticklabels=\empty, yticklabels=\empty
	]
	\foreach \a in {0, 0_1, 0_2, 0_5, 1, 2, 5, 10, 20, 50, 100}{
		\addplot+ [default line] table[header=false, x expr=1/2^\coordindex, y index=3] {data/analysis_B1_a_\a_b_10};
	}
	\foreach \a in {0, 0_1, 0_2, 0_5, 1, 2, 5, 10, 20, 50, 100}{
		\addplot+ [hollow line] table[header=false, x expr=1/2^\coordindex, y index=2] {data/analysis_B1_a_\a_b_10};
	}
	\slopetriangle{.325}{.1}{.075}{3}
	\coordinate (t) at (axis description cs:.985,.03);
	\pgfplotsset{
		after end axis/.append code={%
			\node [anchor=south east] at (t) {$b = 10$};
		}
	}
	\nextgroupplot[%
		xmin=1e-7, xmax=1e0, ymin=1e-16, ymax=1e0, ylabel={$\V{\;\cdot\;}$}, xticklabels=\empty
	]
	\foreach \a in {0, 0_1, 0_2, 0_5, 1, 2, 5, 10, 20, 50, 100}{
		\addplot+ [default line] table[header=false, x expr=1/2^\coordindex, y index=3] {data/analysis_B1_a_\a_b_100};
	}
	\foreach \a in {0, 0_1, 0_2, 0_5, 1, 2, 5, 10, 20, 50, 100}{
		\addplot+ [hollow line] table[header=false, x expr=1/2^\coordindex, y index=2] {data/analysis_B1_a_\a_b_100};
	}
	\slopetriangle{.3}{.1}{.075}{3}
	\slopetriangle{.89}{.1}{.6}{-2}
	\coordinate (t) at (axis description cs:.985,.03);
	\pgfplotsset{
		after end axis/.append code={%
			\node [anchor=south east] at (t) {$b = 100$};
		}
	}
	\nextgroupplot[%
		xmin=1e-7, xmax=1e0, ymin=1e-16, ymax=1e0, yticklabels=\empty, xticklabels=\empty
	]
	\foreach \a in {0, 0_1, 0_2, 0_5, 1, 2, 5, 10, 20, 50, 100}{
		\addplot+ [default line] table[header=false, x expr=1/2^\coordindex, y index=3] {data/analysis_B1_a_\a_b_1000};
	}
	\foreach \a in {0, 0_1, 0_2, 0_5, 1, 2, 5, 10, 20, 50, 100}{
		\addplot+ [hollow line] table[header=false, x expr=1/2^\coordindex, y index=2] {data/analysis_B1_a_\a_b_1000};
	}
	\slopetriangle{.2}{.1}{.075}{3}
	\slopetriangle{.875}{.1}{.425}{-2}
	\coordinate (t) at (axis description cs:.985,.03);
	\pgfplotsset{
		after end axis/.append code={%
			\node [anchor=south east] at (t) {$b = 1\,000$};
		}
	}
	\nextgroupplot[%
		xmin=1e-7, xmax=1e0, ymin=1e-15, ymax=1e-3, ylabel={$\V{\;\cdot\;}$},  xlabel={$\delta t_\level$}
	]
	\foreach \a in {0, 0_1, 0_2, 0_5, 1, 2, 5, 10, 20, 50, 100}{
		\addplot+ [default line] table[header=false, x expr=1/2^\coordindex, y index=3] {data/analysis_B1_a_\a_b_10000};
	}
	\foreach \a in {0, 0_1, 0_2, 0_5, 1, 2, 5, 10, 20, 50, 100}{
		\addplot+ [hollow line] table[header=false, x expr=1/2^\coordindex, y index=2] {data/analysis_B1_a_\a_b_10000};
	}
	\slopetriangle{.225}{.1}{.25}{3}
	\slopetriangle{.875}{.1}{.25}{-2}
	\coordinate (t) at (axis description cs:.985,.03);
	\pgfplotsset{
		after end axis/.append code={%
			\node [anchor=south east] at (t) {$b = 10\,000$};
		}
	}
	\nextgroupplot[%
		xmin=1e-7, xmax=1e0, ymin=1e-15, ymax=1e-3, yticklabels=\empty,  xlabel={$\delta t_\level$}
	]
	\foreach \a in {0, 0_1, 0_2, 0_5, 1, 2, 5, 10, 20, 50, 100}{
		\addplot+ [default line] table[header=false, x expr=1/2^\coordindex, y index=3] {data/analysis_B1_a_\a_b_100000};
	}
	\foreach \a in {0, 0_1, 0_2, 0_5, 1, 2, 5, 10, 20, 50, 100}{
		\addplot+ [hollow line] table[header=false, x expr=1/2^\coordindex, y index=2] {data/analysis_B1_a_\a_b_100000};
	}
	\slopetriangle{.225}{.1}{.4}{3}
	\slopetriangle{.75}{.1}{.15}{-2}
	\coordinate (t) at (axis description cs:.88,.03);
	\pgfplotsset{
		after end axis/.append code={%
			\node [anchor=south east] at (t) {$b = 100\,000$};
		}
	}
	\end{groupplot}
	\path (my plots c1r1.north west|-current bounding box.north) -- coordinate(legendpos) (my plots c2r1.north east|-current bounding box.north);
	\matrix[
		matrix of nodes,
		anchor=south,
		draw=none,
		row sep=1\pgflinewidth,
		inner sep=.2em
	] at ([yshift=1ex]legendpos)
	{
		\ref{V-1-a0} &[0em] $a=0\hphantom{.0}$ &[1em] & \ref{V-1-a0_5} &[0em] $a=0.5$ &[1em] & \ref{V-1-a5} &[0em] $a=5\hphantom{0}$ &[1em] & \ref{V-1-a50} &[0em] $a=50\hphantom{0}$ &[1em]\\
		\ref{V-1-a0_1} &[0em] $a=0.1$ &[1em] & \ref{V-1-a1} &[0em] $a=1\hphantom{.0}$ &[1em] & \ref{V-1-a10} &[0em] $a=10$ &[1em] & \ref{V-1-a100} &[0em] $a=100$ &[1em]\\
		\ref{V-1-a0_2} &[0em] $a=0.2$ &[1em] & \ref{V-1-a2} &[0em] $a=2\hphantom{.0}$ &[1em] & \ref{V-1-a20} &[0em] $a=20$ &[1em] & & & \\
	};
\end{tikzpicture}
	\tikzexternaldisable

\caption{behanviour of the variances $V_\level=\V{\VAqoi_\level}$ (\raisebox{1pt}{\protect\tikz{\protect\draw[default line] (0,0) -- node[inner sep=.8pt, draw, fill=white, circle, pos=.5] {} (.5,0);}}) and $V_{\level, \level-1}=\V{\VAqoi_\level - \VAqoi_{\level-1}}$ (\raisebox{1pt}{\protect\tikz{\protect\draw[default line] (0,0) -- node[inner sep=.8pt, draw, fill, circle, pos=.5] {} (.5,0);}}) for background $\B{1}$ for different values of the parameters $a$ and $b$. Results obtained for $10^6$ particles for $b=1$ and $b=10$, $10^5$ particles for $b=100$ and $b=1\,000$ and $10^4$ particles for $b=10\,000$ and $b=100\,000$.}
\label{fig:V-B1}
\end{figure}

\begin{figure}
\ContinuedFloat
\centering
	\tikzexternalenable
	\tikzsetnextfilename{V-B2}
%
%
\begin{tikzpicture}%
	\begin{groupplot}[%
		group style={
			group size=2 by 3,
			 vertical sep=2em,
			 horizontal sep=2em,
			 group name=my plots
		},
		loglog axis
	]
	
	\nextgroupplot[%
		xmin=1e-7, xmax=1e0, ymin=1e-18, ymax=1e2, ylabel={$\V{\;\cdot\;}$}, xticklabels=\empty
	]
	\foreach \a in {0, 0_1, 0_2, 0_5, 1, 2, 5, 10, 20, 50, 100}{
		\addplot+ [default line] table[header=false, x expr=1/2^\coordindex, y index=3] {data/analysis_B2_a_\a_b_1};
		\label{V-1-a\a}
	}
	\foreach \a in {0, 0_1, 0_2, 0_5, 1, 2, 5, 10, 20, 50, 100}{
		\addplot+ [hollow line] table[header=false, x expr=1/2^\coordindex, y index=2] {data/analysis_B2_a_\a_b_1};
	}
	\slopetriangle{.425}{.1}{.075}{3}
	\coordinate (t) at (axis description cs:.985,.03);
	\pgfplotsset{
		after end axis/.append code={%
			\node [anchor=south east] at (t) {$b = 1$};
		}
	}
	\nextgroupplot[%
		xmin=1e-7, xmax=1e0, ymin=1e-18, ymax=1e2, xticklabels=\empty, yticklabels=\empty
	]
	\foreach \a in {0, 0_1, 0_2, 0_5, 1, 2, 5, 10, 20, 50, 100}{
		\addplot+ [default line] table[header=false, x expr=1/2^\coordindex, y index=3] {data/analysis_B2_a_\a_b_10};
	}
	\foreach \a in {0, 0_1, 0_2, 0_5, 1, 2, 5, 10, 20, 50, 100}{
		\addplot+ [hollow line] table[header=false, x expr=1/2^\coordindex, y index=2] {data/analysis_B2_a_\a_b_10};
	}
	\slopetriangle{.325}{.1}{.075}{3}
	\coordinate (t) at (axis description cs:.985,.03);
	\pgfplotsset{
		after end axis/.append code={%
			\node [anchor=south east] at (t) {$b = 10$};
		}
	}
	\nextgroupplot[%
		xmin=1e-7, xmax=1e0, ymin=1e-16, ymax=1e0, ylabel={$\V{\;\cdot\;}$}, xticklabels=\empty
	]
	\foreach \a in {0, 0_1, 0_2, 0_5, 1, 2, 5, 10, 20, 50, 100}{
		\addplot+ [default line] table[header=false, x expr=1/2^\coordindex, y index=3] {data/analysis_B2_a_\a_b_100};
	}
	\foreach \a in {0, 0_1, 0_2, 0_5, 1, 2, 5, 10, 20, 50, 100}{
		\addplot+ [hollow line] table[header=false, x expr=1/2^\coordindex, y index=2] {data/analysis_B2_a_\a_b_100};
	}
	\slopetriangle{.3}{.1}{.075}{3}
	\slopetriangle{.89}{.1}{.6}{-2}
	\coordinate (t) at (axis description cs:.985,.03);
	\pgfplotsset{
		after end axis/.append code={%
			\node [anchor=south east] at (t) {$b = 100$};
		}
	}
	\nextgroupplot[%
		xmin=1e-7, xmax=1e0, ymin=1e-16, ymax=1e0, yticklabels=\empty, xticklabels=\empty
	]
	\foreach \a in {0, 0_1, 0_2, 0_5, 1, 2, 5, 10, 20, 50, 100}{
		\addplot+ [default line] table[header=false, x expr=1/2^\coordindex, y index=3] {data/analysis_B2_a_\a_b_1000};
	}
	\foreach \a in {0, 0_1, 0_2, 0_5, 1, 2, 5, 10, 20, 50, 100}{
		\addplot+ [hollow line] table[header=false, x expr=1/2^\coordindex, y index=2] {data/analysis_B2_a_\a_b_1000};
	}
	\slopetriangle{.2}{.1}{.075}{3}
	\slopetriangle{.875}{.1}{.425}{-2}
	\coordinate (t) at (axis description cs:.985,.03);
	\pgfplotsset{
		after end axis/.append code={%
			\node [anchor=south east] at (t) {$b = 1\,000$};
		}
	}
	\nextgroupplot[%
		xmin=1e-7, xmax=1e0, ymin=1e-15, ymax=1e-3, ylabel={$\V{\;\cdot\;}$},  xlabel={$\delta t_\level$}
	]
	\foreach \a in {0, 0_1, 0_2, 0_5, 1, 2, 5, 10, 20, 50, 100}{
		\addplot+ [default line] table[header=false, x expr=1/2^\coordindex, y index=3] {data/analysis_B2_a_\a_b_10000};
	}
	\foreach \a in {0, 0_1, 0_2, 0_5, 1, 2, 5, 10, 20, 50, 100}{
		\addplot+ [hollow line] table[header=false, x expr=1/2^\coordindex, y index=2] {data/analysis_B2_a_\a_b_10000};
	}
	\slopetriangle{.225}{.1}{.25}{3}
	\slopetriangle{.875}{.1}{.25}{-2}
	\coordinate (t) at (axis description cs:.985,.03);
	\pgfplotsset{
		after end axis/.append code={%
			\node [anchor=south east] at (t) {$b = 10\,000$};
		}
	}
	\nextgroupplot[%
		xmin=1e-7, xmax=1e0, ymin=1e-15, ymax=1e-3, yticklabels=\empty,  xlabel={$\delta t_\level$}
	]
	\foreach \a in {0, 0_1, 0_2, 0_5, 1, 2, 5, 10, 20, 50, 100}{
		\addplot+ [default line] table[header=false, x expr=1/2^\coordindex, y index=3] {data/analysis_B2_a_\a_b_100000};
	}
	\foreach \a in {0, 0_1, 0_2, 0_5, 1, 2, 5, 10, 20, 50, 100}{
		\addplot+ [hollow line] table[header=false, x expr=1/2^\coordindex, y index=2] {data/analysis_B2_a_\a_b_100000};
	}
	\slopetriangle{.225}{.1}{.4}{3}
	\slopetriangle{.75}{.1}{.15}{-2}
	\coordinate (t) at (axis description cs:.88,.03);
	\pgfplotsset{
		after end axis/.append code={%
			\node [anchor=south east] at (t) {$b = 100\,000$};
		}
	}
	\end{groupplot}
	\path (my plots c1r1.north west|-current bounding box.north) -- coordinate(legendpos) (my plots c2r1.north east|-current bounding box.north);
	\matrix[
		matrix of nodes,
		anchor=south,
		draw=none,
		row sep=1\pgflinewidth,
		inner sep=.2em
	] at ([yshift=1ex]legendpos)
	{
		\ref{V-1-a0} &[0em] $a=0\hphantom{.0}$ &[1em] & \ref{V-1-a0_5} &[0em] $a=0.5$ &[1em] & \ref{V-1-a5} &[0em] $a=5\hphantom{0}$ &[1em] & \ref{V-1-a50} &[0em] $a=50\hphantom{0}$ &[1em]\\
		\ref{V-1-a0_1} &[0em] $a=0.1$ &[1em] & \ref{V-1-a1} &[0em] $a=1\hphantom{.0}$ &[1em] & \ref{V-1-a10} &[0em] $a=10$ &[1em] & \ref{V-1-a100} &[0em] $a=100$ &[1em]\\
		\ref{V-1-a0_2} &[0em] $a=0.2$ &[1em] & \ref{V-1-a2} &[0em] $a=2\hphantom{.0}$ &[1em] & \ref{V-1-a20} &[0em] $a=20$ &[1em] & & & \\
	};
\end{tikzpicture}
	\tikzexternaldisable

\caption{behanviour of the variances $V_\level=\V{\VAqoi_\level}$ (\raisebox{1pt}{\protect\tikz{\protect\draw[default line] (0,0) -- node[inner sep=.8pt, draw, fill=white, circle, pos=.5] {} (.5,0);}}) and $V_{\level, \level-1}=\V{\VAqoi_\level - \VAqoi_{\level-1}}$ (\raisebox{1pt}{\protect\tikz{\protect\draw[default line] (0,0) -- node[inner sep=.8pt, draw, fill, circle, pos=.5] {} (.5,0);}}) for background $\B{2}$ for different values of the parameters $a$ and $b$. Results obtained for $10^6$ particles for $b=1$ and $b=10$, $10^5$ particles for $b=100$ and $b=1\,000$ and $10^4$ particles for $b=10\,000$ and $b=100\,000$.}
\label{fig:V-B2}
\end{figure}

\setlength{\figurewidth}{.45\textwidth}
\setlength{\figureheight}{.45\textwidth}
\begin{figure}
\centering
	\tikzexternalenable
	\tikzsetnextfilename{C-B1}
%
%
\begin{tikzpicture}%
	\begin{groupplot}[%
		group style={
			group size=2 by 3,
			 vertical sep=2em,
			 horizontal sep=2em,
			 group name=my plots
		},
		loglog axis
	]
	\nextgroupplot[%
		xmin=1e-7, xmax=1e0, ymin=1e-6, ymax=1e-1, ylabel={$C_\level$}, xticklabels=\empty
	]
	\foreach \a in {0, 0_1, 0_2, 0_5, 1, 2, 5, 10, 20, 50, 100}{
		\addplot+ [default line] table[header=false, x expr=1/2^\coordindex, y index=4] {data/analysis_B1_a_\a_b_1};
		\label{C-1-a\a}
	}
	\coordinate (t) at (axis description cs:.97,.97);
	\pgfplotsset{
		after end axis/.append code={%
			\node [anchor=north east] at (t) {$b = 1$};
		}
	}
	\nextgroupplot[%
		xmin=1e-7, xmax=1e0, ymin=1e-6, ymax=1e-1, xticklabels=\empty, yticklabels=\empty
	]
	\foreach \a in {0, 0_1, 0_2, 0_5, 1, 2, 5, 10, 20, 50, 100}{
		\addplot+ [default line] table[header=false, x expr=1/2^\coordindex, y index=4] {data/analysis_B1_a_\a_b_10};
	}
	\coordinate (t) at (axis description cs:.97,.97);
	\pgfplotsset{
		after end axis/.append code={%
			\node [anchor=north east] at (t) {$b = 10$};
		}
	}
	\nextgroupplot[%
		xmin=1e-7, xmax=1e0, ymin=1e-6, ymax=1e-1, ylabel={$C_\level$}, xticklabels=\empty
	]
	\foreach \a in {0, 0_1, 0_2, 0_5, 1, 2, 5, 10, 20, 50, 100}{
		\addplot+ [default line] table[header=false, x expr=1/2^\coordindex, y index=4] {data/analysis_B1_a_\a_b_100};
	}
	\slopetriangle{.75}{.1}{.075}{-1}
	\coordinate (t) at (axis description cs:.97,.97);
	\pgfplotsset{
		after end axis/.append code={%
			\node [anchor=north east] at (t) {$b = 100$};
		}
	}
	\nextgroupplot[%
		xmin=1e-7, xmax=1e0, ymin=1e-6, ymax=1e-1, yticklabels=\empty, xticklabels=\empty
	]
	\foreach \a in {0, 0_1, 0_2, 0_5, 1, 2, 5, 10, 20, 50, 100}{
		\addplot+ [default line] table[header=false, x expr=1/2^\coordindex, y index=4] {data/analysis_B1_a_\a_b_1000};
	}
	\slopetriangle{.75}{.1}{.075}{-1}
	\coordinate (t) at (axis description cs:.97,.97);
	\pgfplotsset{
		after end axis/.append code={%
			\node [anchor=north east] at (t) {$b = 1\,000$};
		}
	}
	\nextgroupplot[%
		xmin=1e-7, xmax=1e0, ymin=1e-6, ymax=1e-1, ylabel={$C_\level$},  xlabel={$\delta t_\level$}
	]
	\foreach \a in {0, 0_1, 0_2, 0_5, 1, 2, 5, 10, 20, 50, 100}{
		\addplot+ [default line] table[header=false, x expr=1/2^\coordindex, y index=4] {data/analysis_B1_a_\a_b_10000};
	}
	\slopetriangle{.75}{.1}{.075}{-1}
	\coordinate (t) at (axis description cs:.97,.97);
	\pgfplotsset{
		after end axis/.append code={%
			\node [anchor=north east] at (t) {$b = 10\,000$};
		}
	}
	\nextgroupplot[%
		xmin=1e-7, xmax=1e0, ymin=1e-6, ymax=1e-1, yticklabels=\empty,  xlabel={$\delta t_\level$}
	]
	\foreach \a in {0, 0_1, 0_2, 0_5, 1, 2, 5, 10, 20, 50, 100}{
		\addplot+ [default line] table[header=false, x expr=1/2^\coordindex, y index=4] {data/analysis_B1_a_\a_b_100000};
	}
	\slopetriangle{.75}{.1}{.075}{-1}
	\coordinate (t6) at (axis description cs:.6,.95);
	\coordinate (t) at (axis description cs:.97,.97);
	\pgfplotsset{
		after end axis/.append code={%
			\node [anchor=north east] at (t) {$b = 100\,000$};
		}
	}
	\end{groupplot}
	\path (my plots c1r1.north west|-current bounding box.north) -- coordinate(legendpos) (my plots c2r1.north east|-current bounding box.north);
	\matrix[
		matrix of nodes,
		anchor=south,
		draw=none,
		row sep=1\pgflinewidth,
		inner sep=.2em
	] at ([yshift=1ex]legendpos)
	{
		\ref{C-1-a0} &[0em] $a=0\hphantom{.0}$ &[1em] & \ref{C-1-a0_5} &[0em] $a=0.5$ &[1em] & \ref{C-1-a5} &[0em] $a=5\hphantom{0}$ &[1em] & \ref{C-1-a50} &[0em] $a=50\hphantom{0}$ &[1em]\\
		\ref{C-1-a0_1} &[0em] $a=0.1$ &[1em] & \ref{C-1-a1} &[0em] $a=1\hphantom{.0}$ &[1em] & \ref{C-1-a10} &[0em] $a=10$ &[1em] & \ref{C-1-a100} &[0em] $a=100$ &[1em]\\
		\ref{C-1-a0_2} &[0em] $a=0.2$ &[1em] & \ref{C-1-a2} &[0em] $a=2\hphantom{.0}$ &[1em] & \ref{C-1-a20} &[0em] $a=20$ &[1em] & & & \\
	};
\end{tikzpicture}
	\tikzexternaldisable

\caption{behanviour of the cost $C_\level = \cost{\VAqoi_{\level, \level-1}}$ for background $\B{1}$ for different values of the parameters $a$ and $b$. Results obtained for $10^6$ particles for $b=1$ and $b=10$, $10^5$ particles for $b=100$ and $b=1\,000$ and $10^4$ particles for $b=10\,000$ and $b=100\,000$.}
\label{fig:C-B1}
\end{figure}

\begin{figure}
\ContinuedFloat
\centering
	\tikzexternalenable
	\tikzsetnextfilename{C-B2}
%
%
\begin{tikzpicture}%
	\begin{groupplot}[%
		group style={
			group size=2 by 3,
			 vertical sep=2em,
			 horizontal sep=2em,
			 group name=my plots
		},
		loglog axis
	]
	
	\nextgroupplot[%
		xmin=1e-7, xmax=1e0, ymin=1e-6, ymax=1e-1, ylabel={$C_\level$}, xticklabels=\empty
	]
	\foreach \a in {0, 0_1, 0_2, 0_5, 1, 2, 5, 10, 20, 50, 100}{
		\addplot+ [default line] table[header=false, x expr=1/2^\coordindex, y index=4] {data/analysis_B2_a_\a_b_1};
		\label{C-1-a\a}
	}
	\coordinate (t) at (axis description cs:.97,.97);
	\pgfplotsset{
		after end axis/.append code={%
			\node [anchor=north east] at (t) {$b = 1$};
		}
	}
	\nextgroupplot[%
		xmin=1e-7, xmax=1e0, ymin=1e-6, ymax=1e-1, xticklabels=\empty, yticklabels=\empty
	]
	\foreach \a in {0, 0_1, 0_2, 0_5, 1, 2, 5, 10, 20, 50, 100}{
		\addplot+ [default line] table[header=false, x expr=1/2^\coordindex, y index=4] {data/analysis_B2_a_\a_b_10};
	}
	\coordinate (t) at (axis description cs:.97,.97);
	\pgfplotsset{
		after end axis/.append code={%
			\node [anchor=north east] at (t) {$b = 10$};
		}
	}
	\nextgroupplot[%
		xmin=1e-7, xmax=1e0, ymin=1e-6, ymax=1e-1, ylabel={$C_\level$}, xticklabels=\empty
	]
	\foreach \a in {0, 0_1, 0_2, 0_5, 1, 2, 5, 10, 20, 50, 100}{
		\addplot+ [default line] table[header=false, x expr=1/2^\coordindex, y index=4] {data/analysis_B2_a_\a_b_100};
	}
	\slopetriangle{.75}{.1}{.075}{-1}
	\coordinate (t) at (axis description cs:.97,.97);
	\pgfplotsset{
		after end axis/.append code={%
			\node [anchor=north east] at (t) {$b = 100$};
		}
	}
	\nextgroupplot[%
		xmin=1e-7, xmax=1e0, ymin=1e-6, ymax=1e-1, yticklabels=\empty, xticklabels=\empty
	]
	\foreach \a in {0, 0_1, 0_2, 0_5, 1, 2, 5, 10, 20, 50, 100}{
		\addplot+ [default line] table[header=false, x expr=1/2^\coordindex, y index=4] {data/analysis_B2_a_\a_b_1000};
	}
	\slopetriangle{.75}{.1}{.075}{-1}
	\coordinate (t) at (axis description cs:.97,.97);
	\pgfplotsset{
		after end axis/.append code={%
			\node [anchor=north east] at (t) {$b = 1\,000$};
		}
	}
	\nextgroupplot[%
		xmin=1e-7, xmax=1e0, ymin=1e-6, ymax=1e-1, ylabel={$C_\level$},  xlabel={$\delta t_\level$}
	]
	\foreach \a in {0, 0_1, 0_2, 0_5, 1, 2, 5, 10, 20, 50, 100}{
		\addplot+ [default line] table[header=false, x expr=1/2^\coordindex, y index=4] {data/analysis_B2_a_\a_b_10000};
	}
	\slopetriangle{.75}{.1}{.075}{-1}
	\coordinate (t) at (axis description cs:.97,.97);
	\pgfplotsset{
		after end axis/.append code={%
			\node [anchor=north east] at (t) {$b = 10\,000$};
		}
	}
	\nextgroupplot[%
		xmin=1e-7, xmax=1e0, ymin=1e-6, ymax=1e-1, yticklabels=\empty,  xlabel={$\delta t_\level$}
	]
	\foreach \a in {0, 0_1, 0_2, 0_5, 1, 2, 5, 10, 20, 50, 100}{
		\addplot+ [default line] table[header=false, x expr=1/2^\coordindex, y index=4] {data/analysis_B2_a_\a_b_100000};
	}
	\slopetriangle{.75}{.1}{.075}{-1}
	\coordinate (t6) at (axis description cs:.6,.95);
	\coordinate (t) at (axis description cs:.97,.97);
	\pgfplotsset{
		after end axis/.append code={%
			\node [anchor=north east] at (t) {$b = 100\,000$};
		}
	}
	\end{groupplot}
	\path (my plots c1r1.north west|-current bounding box.north) -- coordinate(legendpos) (my plots c2r1.north east|-current bounding box.north);
	\matrix[
		matrix of nodes,
		anchor=south,
		draw=none,
		row sep=1\pgflinewidth,
		inner sep=.2em
	] at ([yshift=1ex]legendpos)
	{
		\ref{C-1-a0} &[0em] $a=0\hphantom{.0}$ &[1em] & \ref{C-1-a0_5} &[0em] $a=0.5$ &[1em] & \ref{C-1-a5} &[0em] $a=5\hphantom{0}$ &[1em] & \ref{C-1-a50} &[0em] $a=50\hphantom{0}$ &[1em]\\
		\ref{C-1-a0_1} &[0em] $a=0.1$ &[1em] & \ref{C-1-a1} &[0em] $a=1\hphantom{.0}$ &[1em] & \ref{C-1-a10} &[0em] $a=10$ &[1em] & \ref{C-1-a100} &[0em] $a=100$ &[1em]\\
		\ref{C-1-a0_2} &[0em] $a=0.2$ &[1em] & \ref{C-1-a2} &[0em] $a=2\hphantom{.0}$ &[1em] & \ref{C-1-a20} &[0em] $a=20$ &[1em] & & & \\
	};
\end{tikzpicture}
	\tikzexternaldisable

\caption{behanviour of the cost $C_\level = \cost{\VAqoi_{\level, \level-1}}$ for background $\B{2}$ for different values of the parameters $a$ and $b$. Results obtained for $10^6$ particles for $b=1$ and $b=10$, $10^5$ particles for $b=100$ and $b=1\,000$ and $10^4$ particles for $b=10\,000$ and $b=100\,000$.}
\label{fig:C-B2}
\end{figure}

\begin{figure}
\centering%
\begin{tikzpicture}
\draw[very thick, line cap=round, looseness=.45, dashed] (0,0) to[out=50,in=180] (3,3) to[out=0,in=130] (5,1);
\draw[very thick, line cap=round, looseness=.45] (0,0) to[out=50,in=180] (3,3.0005) to[out=0,in=180] (4,2.5) to[out=0, in=210] (5,3);
\draw[very thick, line cap=round, dashed, dash pattern=on 0pt off 4\pgflinewidth] (3,3.5) -- (3,0);
\draw[->, thick, line cap=round] (4.9,1.5) to[out=70,in=290] node[pos=.5, anchor=west] {$a \uparrow$} (4.9,2.75);
\node[anchor=north] at (3,0) {$1/R(x)$};
\node[anchor=west, rotate=90] at (-.35,1.5) {$V_{\level, \level-1}$};
\end{tikzpicture}%
\hspace{4em}
\begin{tikzpicture}
\draw[very thick, line cap=round, rounded corners=6pt] (0,3) -- (3,3) -- (5,0);
\draw[very thick, line cap=round, dashed, dash pattern=on 0pt off 4\pgflinewidth] (3,3.5) -- (3,0);
\node[anchor=north] at (3,0) {$1/R(x)$};
\node[anchor=west, rotate=90] at (-.35,1.5) {$C_{\level, \level-1}$};
\end{tikzpicture}
\caption{behanviour of variance (\emph{left}) and cost (\emph{right}) of the multilevel difference $\Delta \VAqoi_{\level, \level-1}$ with level $\level$. The homogeneous case ($a=0$) is indicated by the dashed line (\raisebox{.15\baselineskip}{\protect\tikz{\protect\draw[very thick, line cap=round, dashed] (0,0) -- (.5,0);}}), and the heterogeneous case ($a\gg0$) is indicated by the full line (\raisebox{.15\baselineskip}{\protect\tikz{\protect\draw[very thick, line cap=round] (0,0) -- (.5,0);}}). }
\label{fig:sketch}
\end{figure}

\begin{figure}
\centering%
\begin{tikzpicture}
\draw[line width=6pt, line cap=round, looseness=.45, red!40] (4.95,1.05) -- (4.95,1.05) (0,0) to[out=50,in=220] (.9,1.1);
\draw[<-, red, thick, line cap=round, shorten <= 10pt, shorten >= 10pt] (.9,1.075) -- (4.95,1.075);
\draw[very thick, line cap=round, looseness=.45, dashed] (0,0) to[out=50,in=180] (3,3) to[out=0,in=130] (5,1);
\draw[very thick, line cap=round, dashed, dash pattern=on 0pt off 4\pgflinewidth] (3,3.5) -- (3,0);
\node[anchor=north] at (3,0) {$1/R(x)$};
\node[anchor=west, rotate=90] at (-.35,1.5) {$V_{\level, \level-1}$};
\end{tikzpicture}%
\hspace{4em}
\begin{tikzpicture}
\draw[line width=6pt, line cap=round, looseness=.45, red!40] (4,2.5) to[out=0, in=210] (5,3) (0,0) to[out=50,in=215] (2.175,2.5);
\draw[<-, red, thick, line cap=round, shorten <= 10pt, shorten >= 10pt] (2.175,2.5) -- (4,2.5);
\draw[very thick, line cap=round, looseness=.45] (0,0) to[out=50,in=180] (3,3.0005) to[out=0,in=180] (4,2.5) to[out=0, in=210] (5,3);
\draw[very thick, line cap=round, dashed, dash pattern=on 0pt off 4\pgflinewidth] (3,3.5) -- (3,0);
\node[anchor=north] at (3,0) {$1/R(x)$};
\node[anchor=west, rotate=90] at (-.35,1.5) {$V_{\level, \level-1}$};
\end{tikzpicture}
\caption{Summary of the level selection strategy from~\secref{subsec:levelselect} for a homogeneous (\emph{left}) and heterogeneous (\emph{right}) background. Levels indicated by \raisebox{-1pt}{\protect\tikz{\protect\draw[line width=6pt, line cap=round, red!40] (0,0) -- (0,0);}} are selected for the heuristic level set $\calL^\star$.}
\label{fig:sketch-with-levels}
\end{figure}

\begin{algorithm}
\begin{algorithmic}[1]
\Statex \textbf{input:} estimates for the variances $V_\level$, $\level=0,1,\ldots,L$ and $V_{\level, \level-1}$, $\level=1,2,\ldots,L$
\Statex \textbf{output:} a set of levels $\calL=\{\level_j\}_{j=1}^J$
\Statex
\Procedure{\texttt{level\_select}}{$V_0, V_1, \ldots, V_L, V_{1, 0}, V_{2, 1}, \ldots, V_{L, L-1}$}

\State $\level \gets 1$
\While{$V_{\level, \level-1} > V_\level$} \Comment{find first level $\level$ where $V_{\level, \level-1} < V_\level$}
	\State $\level \gets \level + 1$
\EndWhile
\State $\calL \gets \{\level-1, \level\}$ \Comment{a minimum of 2 levels is required for MLMC}
\State $V_\text{min} \gets V_{\level, \level-1}$
\For{$\level=0, 1, \ldots, L$}
	\If{$V_{\level, \level-1} < V_\text{min}/2$} \Comment{only add level $\level$ if sufficient variance decay in $V_{\level, \level-1}$}
		\State $\calL \gets \calL \cup \{\level\}$
		\State $V_\text{min} \gets V_{\level, \level-1}$
	\EndIf
\EndFor
\State \Return{$\calL$}
\EndProcedure
\end{algorithmic}
\caption{Heuristic level selection approach.}
\label{alg:level_select}
\end{algorithm}

\subsection{Performance of the ML-KDMC method}\label{subsec:performance}

In this section, we compare the efficiency of the new ML-KDMC scheme to the efficiency of the standard, single-level KDMC scheme in terms of error ($\varepsilon$) versus computational work (wall clock time). We will show numerically that both the single-level and multilevel method have an asymptotic $\varepsilon$-cost complexity of $\order{\varepsilon^{-2}}$, i.e., the expected complexity of an MC-based method, but the constant is significantly reduced for the multilevel scheme. Every experiment consists of an off-line and an on-line part. The off-line part starts by taking $N_\text{warm-up}$ warm-up samples on each level $\level=0,1,\ldots,\tau$, where $\tau$ is such that $\VAtimestep_\tau\approx1/b$, with $b$ the constant parameter in the background collision rate. In our numerical experiments, we used $N_\text{warm-up}=100$. From this set of warm-up samples, we compute the variances $V_\level$, $\level=0,1,\ldots,L$ and $V_{\level, \level-1}$, $\level=1, 2, \ldots, L$, where the values $V_{\level, \level-1}$ with $\level<\tau$ are approximated by the sample variance using the $N_\text{warm-up}$ warm-up samples, and the values $V_{\level, \level-1}$ with $\level>\tau$ are estimated using the asymptotic ratio $V_{\level, \level-1}=\order{\VAtimestep_\level^3}$. These variances are then used as input for the level selection algorithm, i.e.,~\algref{alg:level_select}. It should be stressed that the amount of computational work of this off-line setup phase is orders of magnitude lower than the on-line phase, described below.

In the on-line phase, we run the (single-level) KDMC and multilevel KDMC algorithm repeatedly for a decreasing sequence of tolerances $\varepsilon^{(r)} = 1 / \sqrt{2^r}$ for $r=0, 1, 2, \ldots, R$ imposed on the RMSE. The value $R$ is determined such that the KDMC scheme runs for approximately 10\,000s (wall clock time), and the value $\sqrt{2}$ is such that each simulation takes about twice the amount of work of the previous iteration, assuming the $\varepsilon$-cost complexity of the method scales as $\order{\varepsilon^{-2}}$. The main reason for this $\varepsilon$-adaptive strategy is that it yields more reliable estimates of the bias $|\E{\VAqoi_L-\VAqoi}|$, see, e.g.,~\cite{collier14, robbe17} for details.

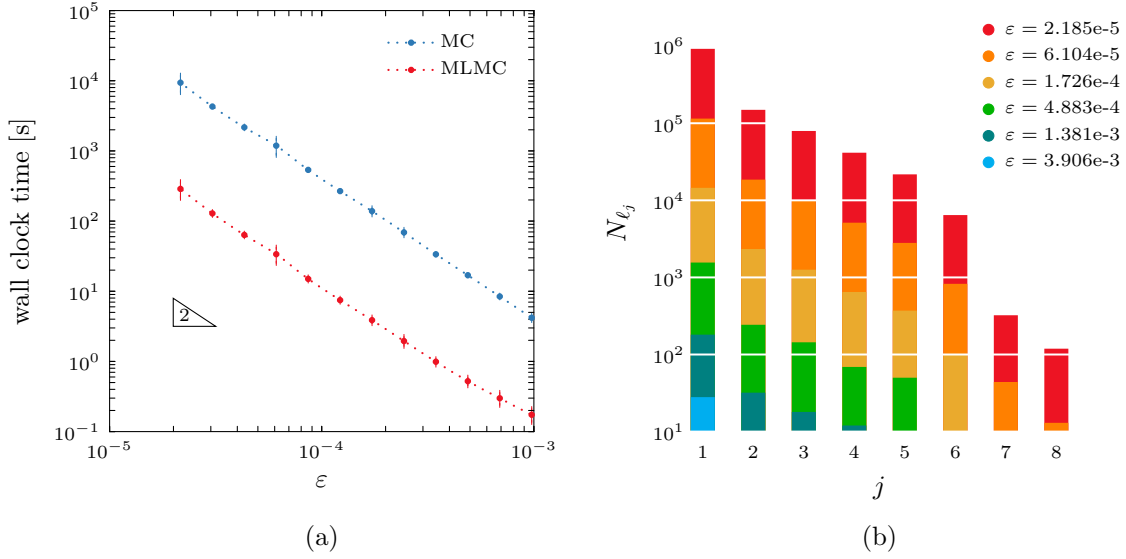
\begin{figure}[p]
	\tikzexternalenable
	\tikzsetnextfilename{N}
\begin{tikzpicture}

	\begin{groupplot}[%
		group style={
			group size=2 by 1,
			 horizontal sep=5em,
			 vertical sep=3em
		}
	]
	\nextgroupplot[%
		default axis,
		xmode=log,
		xlabel={$\varepsilon$},
		xmin=1e-5, xmax=1e-3,
		ymode=log,
		ylabel={wall clock time [s]},
		ymin=1e-1, ymax=1e5,
		legend style={draw=none, font=\scriptsize, at={(0.97,0.97)}, anchor=north east, fill=none, legend cell align=left},
	]
	\addplot[dotted line, lightblue, error bars/.cd, y dir = both, y explicit, error mark options={draw=none}, error bar style={solid}] table[x index=0, y index=1, y error index=2] {figures/complexity.txt};
	\addlegendentry{MC}
	\addplot[dotted line, lightred, error bars/.cd, y dir = both, y explicit, error mark options={draw=none}, error bar style={solid}] table[x index=0, y index=3, y error index=4] {figures/complexity.txt};
	\addlegendentry{MLMC}	
	\slopetriangle{.25}{.1}{.25}{-2}
	\coordinate (b) at (axis description cs:.5,-.2);
	\pgfplotsset{
		after end axis/.append code={%
			\node [anchor=north] at (b) {(a)};
		}
	}
	\nextgroupplot[%
		samples axis, ymode=log, xtick={1, 2, ..., 8}, ymin=10
	]
	\addlegendimage{lightred, bar legend entry}
	\addlegendentry{$\;\varepsilon=2.185\mathrm{\text{e-}}5$}
	\addlegendimage{orange, bar legend entry}
	\addlegendentry{$\;\varepsilon=6.104\mathrm{\text{e-}}5$}
	\addlegendimage{oker, bar legend entry}
	\addlegendentry{$\;\varepsilon=1.726\mathrm{\text{e-}}4$}
	\addlegendimage{lightgreen, bar legend entry}
	\addlegendentry{$\;\varepsilon=4.883\mathrm{\text{e-}}4$}
	\addlegendimage{teal, bar legend entry}
	\addlegendentry{$\;\varepsilon=1.381\mathrm{\text{e-}}3$}
	\addlegendimage{cyan, bar legend entry}
	\addlegendentry{$\;\varepsilon=3.906\mathrm{\text{e-}}3$}
	\addplot[fill=lightred, draw=none] table[x expr=\coordindex+1, y index=5] {figures/samples.txt};
	\addplot[fill=orange, draw=none] table[x expr=\coordindex+1, y index=4] {figures/samples.txt};
	\addplot[fill=oker, draw=none] table[x expr=\coordindex+1, y index=3] {figures/samples.txt};
	\addplot[fill=lightgreen, draw=none] table[x expr=\coordindex+1, y index=2] {figures/samples.txt};
	\addplot[fill=teal, draw=none] table[x expr=\coordindex+1, y index=1] {figures/samples.txt};
	\addplot[fill=cyan, draw=none] table[x expr=\coordindex+1, y index=0] {figures/samples.txt};
	\coordinate (b) at (axis description cs:.5,-.2);
	\pgfplotsset{
		after end axis/.append code={%
			\node [anchor=north] at (b) {(b)};
		}
	}
	\end{groupplot}	
\end{tikzpicture}%
	\tikzexternaldisable

\caption{(a) Cost complexity of ML-KDMC and KDMC for background $\B{1}$ with $a=10$ and $b=1\,000$. (b) Total number of samples on each level $\level_j\in\calL$, $j=1,2,\ldots,8$ for background $\B{1}$ with $a=10$ and $b=1\,000$.}
\label{fig:samples}
\end{figure}

\figref{fig:samples}(a) shows the asymptotic $\varepsilon$-complexity of the ML-KDMC and (single-level) KDMC method for background $\B{1}$ with $a=10$ and $b=1\,000$. Note that both methods indeed follow the asymptotic cost complexity $\order{\varepsilon^{-2}}$. However, the cost of the multilevel method is significantly reduced. Next, in \figref{fig:samples}(b) we show the total number of samples $N_{\level_j}$, $j=1,2,\ldots,J$, in the multilevel method for different values of the tolerance $\varepsilon$ for background $\B{1}$ with $a=10$ and $b=100$. Note that these values are decreasing with increasing $j$. Hence, most samples will be taken with a large time step $\VAtimestep_{\level_1}$, and fewer and fewer samples are required with smaller time step sizes, as claimed in~\secref{subsec:mlmc}. Again, we refer to \url{https://numa.cs.kuleuven.be/research_private/pieterjan.robbe/apmlmc/} for results including all other parameter combinations.

Next, in \tabref{tab:speedup}, we report the algorithmic speedup (computed as the ratio of the amount of computational work expressed in wall clock time) of our ML-KDMC method compared to KDMC for both backgrounds $\B{1}$ and $\B{2}$, and for all background parameter values $a$ and $b$. We note that the MLMC method performs better for larger values of $b$ (i.e., the constant in the background collision rates $R_1(x)$ and $R_2(x)$). The performance decreases slightly with increasing heterogeneity (larger values of $a$). We remark that, in a practical setting, one is interested in large values of the collision rate (large $b$) and moderate to large amounts of heterogeneity (moderate to large values of $a$). We expect our ML-KDMC method to perform very well in these cases.

\setlength\figurewidth{.45\textwidth}
\setlength\figureheight{.45\textwidth}

\begin{table}[p]%
\small%
\centering%
\begin{tabular}{lccccccc}\toprule
& $b = 1$ & $b = 10$ & $b = 100$ & $b = 1\,000$ & $b = 10\,000$ & $b = 100\,000$\\ \midrule
$a=0$ & \cellcolor{yescolor!10}$1.19\pm0.10$ & \cellcolor{yescolor!20}$2.22\pm0.00$ & \cellcolor{yescolor!40}$20.29\pm0.06$ & \cellcolor{yescolor!50}$148.43\pm0.81$ & \cellcolor{yescolor!80}$1\,511.89\pm48.49$ & \cellcolor{yescolor!100}$22\,276.16\pm1\,206.11$ \\$a=0.1$ & \cellcolor{yescolor!10}$1.22\pm0.10$ & \cellcolor{yescolor!20}$2.49\pm0.00$ & \cellcolor{yescolor!40}$20.80\pm2.00$ & \cellcolor{yescolor!60}$162.62\pm0.66$ & \cellcolor{yescolor!80}$1\,196.21\pm20.34$ & \cellcolor{yescolor!100}$14\,076.73\pm\hphantom{0}311.93$ \\$a=0.2$ & \cellcolor{yescolor!10}$1.19\pm0.01$ & \cellcolor{yescolor!20}$2.56\pm0.00$ & \cellcolor{yescolor!40}$21.36\pm1.17$ & \cellcolor{yescolor!60}$162.31\pm0.82$ & \cellcolor{yescolor!80}$1\,203.30\pm13.26$ & \cellcolor{yescolor!100}$11\,931.86\pm\hphantom{0}242.94$ \\$a=0.5$ & \cellcolor{yescolor!10}$1.17\pm0.01$ & \cellcolor{yescolor!20}$2.72\pm0.01$ & \cellcolor{yescolor!40}$21.48\pm0.26$ & \cellcolor{yescolor!50}$118.98\pm0.41$ & \cellcolor{yescolor!70}$1\,000.70\pm\hphantom{0}5.52$ & \cellcolor{yescolor!100}$11\,929.95\pm\hphantom{0}187.13$ \\$a=1$ & \cellcolor{yescolor!10}$1.14\pm0.01$ & \cellcolor{yescolor!20}$2.68\pm0.01$ & \cellcolor{yescolor!30}$17.63\pm0.07$ & \cellcolor{yescolor!50}$120.35\pm0.19$ & \cellcolor{yescolor!70}$1\,000.87\pm\hphantom{0}5.98$ & \cellcolor{yescolor!90}$\hphantom{0}7\,589.55\pm\hphantom{0}168.20$ \\$a=2$ & \cellcolor{yescolor!10}$1.14\pm0.01$ & \cellcolor{yescolor!20}$2.81\pm0.01$ & \cellcolor{yescolor!30}$17.25\pm0.02$ & \cellcolor{yescolor!50}$102.70\pm0.25$ & \cellcolor{yescolor!70}$\hphantom{0}638.25\pm\hphantom{0}1.63$ & \cellcolor{yescolor!100}$\hphantom{0}9\,236.13\pm\hphantom{0}131.17$ \\$a=5$ & \cellcolor{yescolor!10}$1.18\pm0.02$ & \cellcolor{yescolor!20}$2.35\pm0.01$ & \cellcolor{yescolor!30}$15.38\pm0.03$ & \cellcolor{yescolor!50}$\hphantom{0}85.42\pm0.15$ & \cellcolor{yescolor!70}$\hphantom{0}642.67\pm\hphantom{0}1.88$ & \cellcolor{yescolor!100}$\hphantom{0}9\,133.13\pm\hphantom{0}142.74$ \\$a=10$ & \cellcolor{yescolor!10}$1.21\pm0.02$ & \cellcolor{yescolor!20}$1.87\pm0.01$ & \cellcolor{yescolor!30}$14.54\pm0.01$ & \cellcolor{yescolor!50}$\hphantom{0}80.79\pm0.14$ & \cellcolor{yescolor!70}$\hphantom{0}499.99\pm\hphantom{0}1.73$ & \cellcolor{yescolor!90}$\hphantom{0}6\,766.04\pm\hphantom{0}102.58$ \\$a=20$ & \cellcolor{yescolor!10}$1.35\pm0.06$ & \cellcolor{yescolor!20}$2.04\pm0.01$ & \cellcolor{yescolor!30}$17.74\pm0.03$ & \cellcolor{yescolor!50}$\hphantom{0}74.09\pm0.10$ & \cellcolor{yescolor!70}$\hphantom{0}497.18\pm\hphantom{0}2.70$ & \cellcolor{yescolor!90}$\hphantom{0}7\,636.25\pm\hphantom{00}81.69$ \\$a=50$ & \cellcolor{yescolor!10}$1.15\pm0.05$ & \cellcolor{yescolor!20}$2.43\pm0.02$ & \cellcolor{yescolor!30}$\hphantom{0}9.98\pm0.02$ & \cellcolor{yescolor!50}$\hphantom{0}84.74\pm0.15$ & \cellcolor{yescolor!70}$\hphantom{0}457.36\pm\hphantom{0}2.37$ & \cellcolor{yescolor!90}$\hphantom{0}4\,684.73\pm\hphantom{00}77.78$ \\$a=100$ & \cellcolor{yescolor!10}$1.23\pm0.03$ & \cellcolor{yescolor!20}$2.61\pm0.04$ & \cellcolor{yescolor!30}$\hphantom{0}8.69\pm0.05$ & \cellcolor{yescolor!50}$122.60\pm0.66$ & \cellcolor{yescolor!70}$\hphantom{0}452.88\pm\hphantom{0}7.57$ & \cellcolor{yescolor!80}$\hphantom{0}3\,056.32\pm\hphantom{00}27.19$ \\ \bottomrule
\end{tabular}
\caption{Algorithmic speedup of the ML-KDMC method (compared to single-level KDMC) for background $\B{1}$ and various background parameters $a$ and $b$. A green background color (\hspace{-.3pt}\raisebox{-2pt}{\protect\tikz{\protect\node[yes dot] {};}}) is a speedup, a red background color (\hspace{-.3pt}\raisebox{-2pt}{\protect\tikz{\protect\node[no dot] {};}}) is a slowdown. Color saturation indicates the amount of speedup or slowdown.}
\label{tab:speedup}
\end{table}

\begin{table}[p]%
\small%
\ContinuedFloat%
\centering%
\begin{tabular}{lccccccc}\toprule
& $b = 1$ & $b = 10$ & $b = 100$ & $b = 1\,000$ & $b = 10\,000$ & $b = 100\,000$\\ \midrule
$a=0$ & \cellcolor{yescolor!10}$1.19\pm0.09$ & \cellcolor{yescolor!20}$2.22\pm0.01$ & \cellcolor{yescolor!40}$21.29\pm0.04$ & \cellcolor{yescolor!50}$149.07\pm0.34$ & \cellcolor{yescolor!80}$1\,438.15\pm33.72$ & \cellcolor{yescolor!100}$21\,271.88\pm1\,120.02$ \\$a=0.1$ & \cellcolor{yescolor!10}$1.23\pm0.11$ & \cellcolor{yescolor!20}$2.40\pm0.01$ & \cellcolor{yescolor!40}$20.84\pm0.02$ & \cellcolor{yescolor!60}$157.45\pm0.82$ & \cellcolor{yescolor!80}$1\,629.46\pm26.07$ & \cellcolor{yescolor!100}$14\,279.19\pm\hphantom{0}566.27$ \\$a=0.2$ & \cellcolor{yescolor!10}$1.19\pm0.00$ & \cellcolor{yescolor!20}$2.55\pm0.00$ & \cellcolor{yescolor!40}$21.53\pm0.05$ & \cellcolor{yescolor!60}$156.41\pm0.74$ & \cellcolor{yescolor!80}$1\,576.88\pm36.40$ & \cellcolor{yescolor!100}$14\,547.13\pm\hphantom{0}473.81$ \\$a=0.5$ & \cellcolor{yescolor!10}$1.18\pm0.00$ & \cellcolor{yescolor!20}$2.57\pm0.00$ & \cellcolor{yescolor!40}$23.23\pm0.04$ & \cellcolor{yescolor!60}$158.20\pm0.83$ & \cellcolor{yescolor!70}$1\,127.70\pm13.89$ & \cellcolor{yescolor!100}$11\,875.57\pm\hphantom{0}147.67$ \\$a=1$ & \cellcolor{yescolor!10}$1.14\pm0.01$ & \cellcolor{yescolor!20}$2.67\pm0.00$ & \cellcolor{yescolor!40}$21.58\pm0.02$ & \cellcolor{yescolor!50}$118.54\pm0.28$ & \cellcolor{yescolor!70}$\hphantom{0}991.45\pm\hphantom{0}7.17$ & \cellcolor{yescolor!90}$\hphantom{0}7\,446.33\pm\hphantom{0}135.31$ \\$a=2$ & \cellcolor{yescolor!10}$1.16\pm0.02$ & \cellcolor{yescolor!20}$2.66\pm0.00$ & \cellcolor{yescolor!40}$19.91\pm0.05$ & \cellcolor{yescolor!50}$111.60\pm0.36$ & \cellcolor{yescolor!70}$\hphantom{0}986.87\pm\hphantom{0}5.98$ & \cellcolor{yescolor!100}$\hphantom{0}9\,181.21\pm\hphantom{0}142.87$ \\$a=5$ & \cellcolor{yescolor!10}$1.41\pm0.16$ & \cellcolor{yescolor!20}$2.67\pm0.00$ & \cellcolor{yescolor!30}$17.10\pm0.02$ & \cellcolor{yescolor!50}$\hphantom{0}92.69\pm0.23$ & \cellcolor{yescolor!70}$\hphantom{0}632.20\pm\hphantom{0}1.71$ & \cellcolor{yescolor!100}$\hphantom{0}9\,352.35\pm\hphantom{0}160.51$ \\$a=10$ & \cellcolor{yescolor!10}$1.17\pm0.01$ & \cellcolor{yescolor!20}$2.40\pm0.00$ & \cellcolor{yescolor!30}$14.43\pm0.01$ & \cellcolor{yescolor!50}$\hphantom{0}66.72\pm0.07$ & \cellcolor{yescolor!70}$\hphantom{0}635.49\pm\hphantom{0}2.81$ & \cellcolor{yescolor!100}$\hphantom{0}9\,250.00\pm\hphantom{0}161.89$ \\$a=20$ & \cellcolor{nocolor!10}$1.00\pm0.04$ & \cellcolor{yescolor!20}$2.41\pm0.02$ & \cellcolor{yescolor!30}$14.24\pm0.01$ & \cellcolor{yescolor!50}$\hphantom{0}80.74\pm0.11$ & \cellcolor{yescolor!70}$\hphantom{0}505.88\pm\hphantom{0}1.49$ & \cellcolor{yescolor!100}$\hphantom{0}8\,605.50\pm\hphantom{0}133.31$ \\$a=50$ & \cellcolor{yescolor!10}$1.25\pm0.00$ & \cellcolor{yescolor!20}$1.82\pm0.05$ & \cellcolor{yescolor!30}$15.26\pm0.03$ & \cellcolor{yescolor!50}$\hphantom{0}70.84\pm0.13$ & \cellcolor{yescolor!70}$\hphantom{0}484.00\pm\hphantom{0}2.07$ & \cellcolor{yescolor!90}$\hphantom{0}6\,159.53\pm\hphantom{0}136.14$ \\$a=100$ & \cellcolor{nocolor!10}$1.00\pm0.00$ & \cellcolor{yescolor!10}$1.68\pm0.01$ & \cellcolor{yescolor!30}$18.99\pm0.07$ & \cellcolor{yescolor!50}$\hphantom{0}66.06\pm0.13$ & \cellcolor{yescolor!70}$\hphantom{0}450.02\pm\hphantom{0}4.35$ & \cellcolor{yescolor!90}$\hphantom{0}5\,019.00\pm\hphantom{0}107.72$ \\ \bottomrule
\end{tabular}
\caption{Algorithmic speedup of the ML-KDMC method (compared to single-level KDMC) for background $\B{2}$ and various background parameters $a$ and $b$. A green background color (\hspace{-.3pt}\raisebox{-2pt}{\protect\tikz{\protect\node[yes dot] {};}}) is a speedup, a red background color (\hspace{-.3pt}\raisebox{-2pt}{\protect\tikz{\protect\node[no dot] {};}}) is a slowdown. Color saturation indicates the amount of speedup or slowdown.}
\end{table}

\section{Conclusion and future work}

In this paper, we introduced a multilevel extension of the Kinetic-Diffusion Monte Carlo (KDMC) scheme for solving the Boltzmann-BGK equation proposed in~\cite{mortier20}, called Multilevel Kinetic-Diffusion Monte Carlo (ML-KDMC). Crucial in our algorithm is the new and improved recipe for correlated sampling of a particle trajectory with different time step sizes. We show that this correlated sampling can be achieved using a mapping and aggregation strategy for the random numbers used in the particle path simulation. We also discussed specific challenges in applying the multilevel sampling strategy to the KDMC scheme. Notably, because of the nonmonotone behanviour of the variance and cost of the multilevel differences as a function of the time step size, the selection of the appropriate hierarchy of larger time step sizes becomes a nontrivial problem. We introduced a heuristic method for this level selection problem, that avoids solving a combinatorial optimization problem involving all possible combinations of larger time step sizes. We illustrate numerically that our new ML-KDMC scheme with optimal level hierarchy outperforms the classic, single-level KDMC scheme in terms of error versus computational cost by several orders of magnitude, for a wide variety of background collision rate parameter combinations. The implementation of the (ML-)KDMC scheme in actual nuclear fusion plasma simulation codes, such as EIRENE~\cite{reiter05}, is the topic of currently ongoing research.

Finally, we remark that our ML-KDMC simulation scheme is not limited to the kinetic simulation of the Boltzmann-BGK equation in nuclear fusion. Applications in, amongst others, rarefied gases, see~\cite{pareschi99}, and radiation transport, see~\cite{fleck71}, can also benefit from our hybrid (ML-)KDMC method.

\bibliographystyle{acm}
\bibliography{references}

\end{document}